\documentclass[12pt, reqno, oneside]{amsart}

\usepackage{amsmath, amssymb, amscd, amsthm, comment, amsxtra}
\usepackage{graphicx}
\usepackage{subfigure}
\usepackage{enumerate}
\usepackage[colorlinks=true,urlcolor=blue]{hyperref}
\usepackage{euscript}
\usepackage[format=plain,labelfont=bf,up,width=.99\textwidth]{caption}
\usepackage[toc,page]{appendix}
\usepackage{mathabx}
\usepackage{xcolor}

\usepackage{tikz}

\newcounter{reminder}

\newcommand{\hide}[1]{}

\theoremstyle{plain}
\newtheorem{thm}{Theorem}[section]
\newtheorem*{mainthm}{Main Theorem}

\newtheorem{lem}[thm]{Lemma}
\newtheorem{prop}[thm]{Proposition}
\newtheorem{cor}[thm]{Corollary}

\theoremstyle{definition}
\newtheorem{defn}[thm]{Definition}

\newtheorem{rem}[thm]{Remark}

\renewcommand{\arg}{\operatorname{arg}}

\renewcommand{\mod}{\operatorname{mod}}

\newcommand{\Crit}{\operatorname{Crit}}

\newcommand{\crit}{\operatorname{crit}}

\newcommand{\gen}{\operatorname{gen}}

\newcommand{\Int}{\operatorname{int}}
\newcommand{\rt}{\operatorname{root}}

\numberwithin{equation}{section}
\newcommand{\thmref}[1]{Theorem~\ref{#1}}
\newcommand{\propref}[1]{Proposition~\ref{#1}}
\newcommand{\lemref}[1]{Lemma~\ref{#1}}

\newcommand{\bbN}{\mathbb N}
\newcommand{\bbZ}{\mathbb Z}
\newcommand{\bbQ}{\mathbb Q}
\newcommand{\bbR}{\mathbb R}
\newcommand{\bbC}{\mathbb C}
\newcommand{\bbD}{\mathbb D}

\newcommand{\Bp}{\mathbf{P}}

\newcommand{\Cc}{\widehat{\mathbb C}}

\newcommand{\fs}{\mathfrak{s}}
\newcommand{\ft}{\mathfrak{t}}
\newcommand{\fm}{\mathfrak{m}}

\newcommand{\tf}{\tilde{f}}

\newcommand{\Pp}{\mathcal P}

\newcommand{\tQ}{\tilde{\mathcal Q}}

\newcommand{\cA}{\mathcal A}

\newcommand{\cE}{\mathcal E}

\newcommand{\cJ}{\mathcal J}

\newcommand{\cL}{\mathcal L}

\newcommand{\cO}{\mathcal O}

\newcommand{\cQ}{\mathcal Q}
\newcommand{\cR}{\mathcal R}
\newcommand{\cS}{\mathcal S}
\newcommand{\cT}{\mathcal T}

\newcommand{\cV}{\mathcal V}

\newcommand{\U}{\mathcal U}
\newcommand{\V}{\mathcal V}
\newcommand{\tV}{\widetilde{\mathcal V}}
\newcommand{\C}{\mathbb C}

\newcommand{\Y}{\mathcal Y}
\newcommand{\X}{\mathcal X}

\newcommand{\Kwi}{K_{\text{well-inside}}}

\newcommand{\Koc}{K_{\text{off-crit}}}

\newcommand{\sm}{\setminus}

\renewcommand{\tilde}{\widetilde}

\DeclareMathOperator{\fib}{fib}
\DeclareMathOperator{\orb}{orb}
\DeclareMathOperator{\Comp}{Comp}
\DeclareMathOperator{\PC}{PC}
\DeclareMathOperator{\modulus}{mod}
\DeclareMathOperator{\cl}{cl}
\DeclareMathOperator{\area}{area}
\DeclareMathOperator{\inter}{int}

\newcommand{\cW}{\mathcal W}
\newcommand{\cX}{\mathcal X}
\newcommand{\cY}{\mathcal Y}
\newcommand{\cZ}{\mathcal Z}

\newcommand{\bfR}{\mathbf R}

\newcommand{\bff}{\mathbf f}

\newcommand{\tif}{\tilde f}

\newcommand{\ticE}{\tilde \cE}

\newcommand{\tiV}{\tilde V}

\newcommand{\tieta}{\tilde \eta}

\DeclareMathOperator{\Dom}{\mathcal R} 
\DeclareMathOperator{\DomE}{\mathcal L}		
\DeclareMathOperator{\DomL}{\hat{\mathcal L}}	

\newcommand{\tic}{\tilde c}
\newcommand{\tih}{\tilde h}
\newcommand{\tiT}{{\tilde T}}
\newcommand{\tiv}{{\tilde v}}

\newcommand{\tiu}{\tilde u}
\newcommand{\tie}{\tilde e}

\newcommand{\tipsi}{\tilde \psi}

\newcommand{\ticT}{\tilde \cT}
\newcommand{\ticV}{\tilde \cV}

\newcommand{\ticO}{\tilde \cO}
\newcommand{\tiW}{\tilde W}

\newcommand{\adr}{\operatorname{adr}}
\newcommand{\het}{\operatorname{height}}
\newcommand{\jt}{\operatorname{jt}}

\newcommand{\matsp}[1]{\hspace{5mm} \text{#1} \hspace{5mm}}

\title[Rigidity of Bounded-type Siegel Polynomials]{Rigidity of Bounded-type Siegel Polynomials}
\author{Kostiantyn Drach}
\author{Jonguk Yang}

\address{Universitat de Barcelona, Gran Via 585, 08007 Barcelona, Spain}
\address{Centre de Recerca Matem\`atica, Edifici C, Carrer de l'Albareda, 08193 Bellaterra, Barcelona, Spain}
\email{kostiantyn.drach@ub.edu}
\address{Korea Institute for Advanced Study, 85 Hoegi-ro, Dongdaemun-gu, Seoul 02455, Republic of Korea}
\email{jongukyang@kias.re.kr}

\thanks{The first author is partially supported by Departament de Recerca i Universitats de la Generalitat de Catalunya (2021 SGR 00697), Agencia Estatal de Investigaci\'on grants PID2023-147252NB-I00, CNS2025-166633, and Maria de Maeztu Excellence Grant CEX2020-001084-M. The second author was partially supported by June E Huh Center for Mathematical Challenges at the Korea Institute of Advanced Studies.}

\begin{document}

\begin{abstract} 
We establish rigidity for a class of higher-degree complex polynomials with
irrationally indifferent dynamics. Specifically, we consider non-renormalizable
(in the sense of Douady and Hubbard) polynomials of degree \(d\geq 2\) with a
Siegel disk whose rotation number is of bounded type. We call such maps
\emph{atomic Siegel polynomials of bounded type}. Our main results are:

\smallskip
\begin{enumerate}
\item[(A)]
The Julia set of every atomic Siegel polynomial of bounded type is locally
connected.

\item[(B)]
Every atomic Siegel polynomial of bounded type is quasiconformally rigid;
equivalently, its Julia set supports no invariant line fields.

\item[(C)]
Any two combinatorially equivalent atomic Siegel polynomials of bounded type
are affinely conjugate.
\end{enumerate}
\smallskip

In particular, (C) proves the \emph{Combinatorial Rigidity Conjecture} for atomic Siegel polynomials of bounded type in arbitrary degree. This extends the higher-degree rigidity theory of Avila--Kahn--Lyubich--Shen \cite{AKLS} and Kozlovski--van Strien \cite{KvS} to the setting of irrationally indifferent dynamics, a setting not previously covered by Yoccoz-type rigidity results.
\end{abstract}

\maketitle

\setcounter{tocdepth}{1}
\tableofcontents

\section{Introduction}\label{sec:intro}

A central theme in one-dimensional complex dynamics is the study of \emph{rigidity phenomena}, that is, understanding to what extent the combinatorial or topological data of a dynamical system determine its analytic structure. This perspective has driven substantial progress in the field, beginning with the seminal work of Douady and Hubbard \cite{DoHu}. They laid out a framework for describing the combinatorial structure of polynomial dynamics through the landing patterns of external rays. Moreover, their theory of polynomial-like mappings and renormalization made it possible to isolate recurrent dynamics at smaller scales and study them separately as independent systems. In this way, they provided a mechanism for decomposing complicated dynamics into simpler constituent pieces.

A decisive step toward rigidity was later taken by Yoccoz \cite{HY}. He pioneered the use of \emph{puzzles}: dynamically defined partitions obtained from external rays and equipotentials. Successive pullbacks of these partitions produce nested neighborhoods on points in the Julia set that record the combinatorial itinerary of an orbit at increasingly fine scales. For non-renormalizable quadratic polynomials with all periodic points repelling, Yoccoz proved that every nested sequence of puzzle pieces shrinks to a point. This shrinking theorem provides the geometric control needed to promote combinatorial equivalence to rigidity, and has since become a foundational tool in complex dynamics.

For Yoccoz's argument, it is essential that the polynomial is quadratic. This ensures that all losses of modulus under pullback are governed by a single simple critical point. For polynomials of higher degree, higher-order branching and the possible interaction of several critical orbits prevent a direct extension of the quadratic method. Generalizing the theory to the higher degree case therefore required more powerful analytic and combinatorial tools, including the Covering Lemma of Kahn and Lyubich \cite{KaLy1,KaLy2}, refined puzzle constructions (such as the enhanced nest of Kozlovski--Shen--van Strien \cite{KSS}), and the theory of generalized renormalizations. Through the work of Lyubich \cite{Ly}, Avila--Kahn--Lyubich--Shen \cite{AKLS}, Kozlovski--Shen--van Strien \cite{KSS}, Kozlovski--van Strien \cite{KvS}, and others, Yoccoz's rigidity result was ultimately extended to non-renormalizable polynomials of arbitrary degree.\footnote{For a recent survey of developments toward rigidity in the quadratic family, see Dudko \cite{Dima}.} Indifferent periodic orbits, however, remained outside the scope of this theory.

The goal of this paper is to extend Yoccoz's rigidity result to polynomials with irrationally indifferent periodic orbits whose rotation numbers are of bounded type. Let $f \colon \C \to \C$ be a polynomial of degree $\deg(f) \geq 2$ with a fixed point at $0$. Denote its Julia set and Fatou set by $J(f)$ and $F(f) := \C \setminus J(f)$ respectively. We say that $0$ is {\it irrationally indifferent} if
$$
f'(0) = e^{i 2\pi \theta}
\matsp{for some}
\theta \in (\bbR \setminus \bbQ)/\bbZ.
$$
The angle $\theta$ is referred to as the {\it rotation number}. We say that $0$ is a {\it Siegel point} if $f$ is locally linearizable near $0$ (or equivalently, if $0 \in F(f)$). Otherwise, $0$ is called a {\it Cremer point}. In the former case, the linearizing conjugacy has a maximal extension to a conformal map from the Fatou component $\Delta_f \subset F(f)$ containing $0$ to the standard disk $\bbD$. The set $\Delta_f$ is referred to as a {\it Siegel disk}. If $f$ has a Siegel fixed point at $0$, it is referred to as a {\it Siegel polynomial (with rotation number $\theta$)}. We say that the rotation number \(\theta\) is of \emph{bounded type} if the entries in its continued-fraction expansion are uniformly bounded. It is known that an irrationally indifferent periodic point is Siegel whenever its rotation number is of bounded type \cite{Si}.

The problem of rigidity becomes particularly difficult in the presence of irrationally indifferent periodic orbits. Near a neutral orbit, puzzle pieces may undergo arbitrarily many iterations
without acquiring definite expansion, and the usual hyperbolic mechanisms
for proving their shrinking do not apply. At the same time, the boundary of a Siegel disk necessarily interacts with recurrent critical dynamics. In this highly non-linear situation, it is a priori possible for the geometry to become extremely distorted, leading to a failure of rigidity.

A major early breakthrough in the study of Siegel disks was due to the combined efforts of Douady, Ghys, Herman, Shishikura, and \'Swi\c atek. They discovered that a quadratic polynomial with a Siegel disk of bounded type rotation number can be modeled by a Blaschke product via quasiconformal surgery. This led to the celebrated result that the boundaries of such Siegel disks are quasi-circles containing the critical point \cite{Do,He2}. By analyzing the surgery map defined on a certain space of some degree-5 Blaschke products, Zakeri proved that the analog of the above results is also true for cubic polynomials \cite{Za}. For polynomials of arbitrary degree, Shishikura announced that bounded-type Siegel boundaries are quasi-circles and that each periodic cycle of such
boundaries contains a critical point. G. Zhang later extended this result to
all rational maps \cite{Zh}. Building on these ideas, the second author constructed Blaschke product models for bounded-type Siegel polynomials of arbitrary degree \cite{Y}.

In the phase space of a Blaschke product model, the Siegel boundary is straightened to the unit circle. The resulting symmetry greatly facilitates geometric control of the system. In particular, it enables the use of the renormalization theory and a priori bounds developed for critical circle maps by numerous authors, including Herman \cite{He2}, Yampolsky \cite{Ya1}, de Faria--de Melo \cite{dFdM}, and Estevez--Smania--Yampolsky \cite{EsSmYa}.

Working with the Blaschke product models of quadratic Siegel polynomials, Petersen proved that their Julia sets are locally connected and have zero area \cite{Pe}; see also Yampolsky's alternative approach using complex a priori bounds \cite{Ya1}. These results established rigidity in the quadratic case for bounded-type rotation numbers.

In Petersen's work, it is crucial that the polynomial is quadratic, since its unique critical point is then forced to lie on the Siegel boundary. For a polynomial \(f\) of higher degree, however, no such restriction applies to its additional critical points. Nonetheless, the orbits of these  \emph{free critical points} can still recur arbitrarily close to \(\partial\Delta_f\), and the resulting distortion is not controlled by the critical circle-map theory. A priori, this may give rise to a \emph{Siegel hedgehog}: a one-sided hairy Siegel disk with a free critical point at the tip of one of its ``hairs.''

This local obstruction is ruled out by the second author in \cite{Y}. More precisely, it was shown that a bounded-type Siegel disk $\Delta_f$ of a polynomial $f$ is not contained in any larger invariant Siegel continuum, and local connectivity holds at the boundary $\partial \Delta_f$. Thus, the topology of the Julia set is rigid at \(\partial\Delta_f\), even in the presence of free critical orbits.

We now turn from this local result to the problem of global rigidity for higher-degree Siegel polynomials. We restrict our attention to Siegel polynomials that admit no nontrivial polynomial-like renormalization with connected filled Julia set. We call such polynomials \emph{atomic}. This assumption isolates the dynamically indecomposable component containing the Siegel dynamics, thereby capturing the essential interaction between the
neutral dynamics and the free critical orbits. The presence of the neutral orbit prevents a direct application of the existing higher-degree machinery developed for Yoccoz-type results and requires substantial modification. In the cubic case, where there is only one free critical point, this problem was resolved by the second author and R.~Zhang \cite{YZh}. The purpose of the present paper is to establish global rigidity in arbitrary degree.

Our main result is the following.

\begin{mainthm}[Rigidity of Siegel polynomials of bounded type]
\label{Thm:A}
Let $f \colon \C \to \C$ be an atomic Siegel polynomial of bounded type and degree $d \ge 2$. Then: 

\begin{enumerate}[\rm (A)]
\item
\label{it:main:1}
each point in the Julia set \(J_f\) is dynamically rigid (i.e., distinct points have distinct combinatorial itineraries); in particular, $J_f$ is locally connected.
\item
\label{it:main:2}
$J_f$ is quasiconformally rigid (i.e., it does not admit non-trivial invariant line fields).
\item
\label{it:main:3}
If $\tilde f$ is another atomic Siegel polynomial of bounded-type that is combinatorially equivalent to $f$, then $f$ and $\tilde f$ are affinely conjugate.
\end{enumerate}
\end{mainthm}

We will define all the concepts in the statement of the Main Theorem in the follow-up sections. In particular, the combinatorial equivalence is defined in Section~\ref{sec:topmodel} (see Definition~\ref{Def:CombEq}). Our combinatorial data includes the rotation number, addresses of critical points in the ``Siegel tree", as well as conformal positions of the capture critical orbits that eventually land on the Siegel disk.

Main Theorem~\eqref{it:main:3} confirms the \emph{Combinatorial Rigidity Conjecture} in the above-mentioned class of polynomials. This conjecture asserts that for polynomials under an \emph{appropriate} definition of the combinatorial equivalence, combinatorial rigidity holds. We remark that our definition includes the condition that $f$ and $\tilde f$ are \emph{Fatou normalized} and the positions of critical points on the boundary of the Siegel disk are combinatorially the same; see Section~\ref{sec:fatoutrees} for details. Since we assume that $f, \tilde f$ are non-renormalizable, their rational laminations, i.e., the combinatorial pattern of landing external rays, are empty. In our context, Fatou normalization means that each critical point $c$ in a preimage of a Siegel disk $\Delta_f$ lands at $0$.

\begin{rem}
The Main Theorem extends readily to polynomials that are at most finitely renormalizable in the sense of Douady--Hubbard. This includes, in particular, polynomials for which several bounded-type Siegel disks and attracting basins coexist.

A less immediate extension would be to allow rationally indifferent periodic points. The difficulty is that the external rays used to separate the dynamical components containing such cycles need not co-land at repelling periodic points. Consequently, the induced renormalizations need not be polynomial-like in the strict sense of Douady--Hubbard. Nevertheless, we expect rigidity to continue to hold and the present proof to require only minor modifications.
\end{rem}

\subsection{Outline of the proof}

The proof uses the following strategy, which was successful in a number of cases (see \cite[Section 2]{CDKvS}), including the one in the current paper:

\begin{enumerate}
\item
\emph{Markov partition}: In a given family $\mathcal F$ of mappings and a map $f \in \mathcal F$ in this family, construct a puzzle partition (i.e., a Yoccoz- or Markov-type partition) around points in the Julia set $J_f$ of $f$. In the context of the current paper, the partition is given in terms of Siegel bubble rays extended by external rays, see Section~\ref{sec:puzzle}.

\item
\emph{Rigidity on the Julia boundary of puzzle pieces}: Establish (dynamical) rigidity on the boundary of each puzzle piece, i.e., for a puzzle piece $P$, show that each point in $\partial P \cap J_f$ can be surrounded by a union of arbitrary small puzzle pieces. In our case, this result was proven in \cite{Y} (see Section~\ref{sec:puzzle}).

\item
\emph{Extract a dynamically natural complex box mapping}: Construct a dynamically natural box mapping (in the sense of Section~\ref{sec:box}) around the critical points of $f$ that do not land on the puzzle boundary. This mapping controls the dynamics of most critical points, except for those whose orbits accumulate only at the puzzle boundary. In our paper, the construction of such a mapping is carried out in Section~\ref{sec:DNBM}.

\item
\emph{Local connectivity and quasiconformal rigidity of the Julia set}: Use the rigidity result for box mappings (i.e., puzzle pieces shrink to points, no invariant line fields) to infer the corresponding results for $f \in \mathcal F$. In our case, this is done in Section~\ref{sec:proofAB}, where we prove the items \eqref{it:main:1} and \eqref{it:main:2} of the Main Theorem.

\item 
\emph{Topological model of the Julia set}: Construct topological models of the dynamics for maps $f \in \mathcal F$ and prove that they are determined uniquely by combinatorial data. This is done in Section~\ref{sec:combtrees} and \ref{sec:fatoutrees}.

\item 
\emph{Topological rigidity of combinatorially equivalent maps}: Prove that the filled Julia sets of maps in $\mathcal F$ are homeomorphic to the topological models constructed in the previous step. It follows that two maps within the same combinatorial class must be topologically conjugate. This is carried out in Section~\ref{sec:topmodel}.

\item
\emph{Conformal rigidity of topologically equivalent maps}: Use quasiconformal rigidity of the corresponding complex box mappings for the topologically conjugate $f, \tilde f$ to conclude that $f, \tilde f$ are quasiconformally conjugate, and hence, because of quasiconformal rigidity of $J_f, J_{\tilde f}$ and Fatou normalization, $f, \tilde f$ are conformally (affine) conjugate. 

In our paper, this step is done in Section~\ref{sec:top2qc}. This concludes the proof of Main Theorem~\eqref{it:main:3}.   
\end{enumerate}

\section{Some preliminaries}

Let $f \colon \Cc \to \Cc$ be a polynomial of degree $d \geq 2$. Suppose that $0$ is an irrationally indifferent fixed point that is the dynamical center of a Siegel disk $\Delta_f$ with the rotation number $\rho = \rho(f) \in (\bbR \setminus \bbQ)/\bbZ$ of bounded type. We refer to such mappings as \emph{Siegel polynomials}.

Denote the attracting basin of infinity, the filled Julia set and the Julia set of $f$ by $\cA^\infty_f$,
$$
K_f := \bbC \setminus \cA^\infty_f
\matsp{and}
J_f := \partial K_f = \partial \cA^\infty_f.
$$

Let $\phi^\infty_f : \cA^\infty_f \to \bbC \setminus \overline{\bbD}$ be the external B\"ottcher uniformization which conjugates $f|_{\cA^\infty_f}$ with $z \mapsto z^d$. The {\it external ray} with {\it external angle} $t\in \bbR/\bbZ$ is given by
$$
\cR^\infty_t := \{z \in \cA^\infty_f \; | \; \arg(\phi^\infty_f(z)) = t\}.
$$
The external ray $\cR^\infty_t$ is said to be {\it rational} if $t \in \bbQ$, and {\it irrational} if otherwise. The accumulation set of $\cR^\infty_t$ is denoted $\omega(\cR^\infty_t) \subset \cJ_f$. We say that $\cR^\infty_t$ {\it lands} at $x$ if $\omega(\cR^\infty_t) = \{x\}$. It is well-known that if $t$ is periodic, then $\cR^\infty_t$ lands at either a repelling or a parabolic periodic point.

We denote by $\Crit(f)$ the set of critical points of $f$ and say that a Siegel polynomial $f$ is \emph{polynomial-like renormalizable around a critical point $c \in \Crit(f)$} if there exist a pair of topological disks $U \Subset V \subset \bbC$ and an integer $n \in \bbN$ such that $c \in U$ and $f^n \colon U  \to V$ is a branched covering such that its \emph{non-escaping set} (or \emph{filled Julia set}) $\{z \in U : f^{nk}(z) \in U, \ k \geqslant 1\}$  is connected and contains $c$. We exclude trivial renormalizations, when the non-escaping set is equal to $K_f$. This is a classical definition due to Douady and Hubbard, and the map $f^n \colon U \to V$ is called a \emph{polynomial-like restriction of $f$} (around $c$). The map $f$ is \emph{non-polynomial-like renormalizable}, or more simply \emph{atomic}, if $f$ admits no nontrivial polynomial-like restriction around any of its critical points.

Finally, the polynomial $f$ is at \emph{most finitely many times polynomial-like renormalizable} if $f$ has at most finitely many polynomial-like restrictions around some critical point in the sense given above (with the non-escaping set getting strictly smaller with each further renormalization).

Observe that if $f$ is an atomic Siegel polynomial, then the Julia set of $f$ is connected. This is because otherwise one might extract a polynomial-like restriction containing $\Delta_f$.

\subsection{Siegel Puzzle Partition}\label{sec:puzzle}

In this section, we provide a puzzle partition suitable for Siegel polynomials. The results of this subsection, unless otherwise stated, are proven in \cite{Y}.

\begin{thm}[Shishikura, Zhang]\label{siegel bound quasi}
The Siegel boundary $\partial \Delta_f$ is a quasi-circle that contains at least one critical point $c_0$ of $f$. \qed
\end{thm}

\begin{cor}
There exists a quasi-symmetric map $h \colon (\partial \Delta_f, c_0) \to (\partial \bbD, 1)$ such that for $g := f|_{\partial \Delta_f}$, we have
$$
h \circ g \circ h^{-1}(z) = e^{2\pi \rho i}z
\matsp{for}
z \in \partial \Delta_f.
$$
\qed
\end{cor}

For $s \in \bbR/\bbZ$, let
$$
\xi_s := h^{-1}(e^{2\pi s i}).
$$
For $k \in \bbZ$, denote
$$
c_k := g^k(c_0) = \xi_{k\rho}.
$$
Without loss of generality, we may assume that $c_k$ is not a critical point for $k \geq 1$.

A connected component $B$ of $f^{-n}(\Delta_f)$ for $n \geq 0$ is called a {\it bubble}. Its {\it generation} $\gen(B)$ is the smallest number $k$ such that $f^k(B) = \Delta_f$. We assume that $\Delta_f$ is the unique bubble of generation $0$. For $\gen(B) \ge 1$, a {\it root} of $B$ is a point in $f^{-\gen(B)+1}(c_0) \cap \partial B$. Note that a bubble can have several roots.

Let $\{B_i\}_{i=0}^\infty$ be a sequence of bubbles such that
\begin{itemize}
\item $B_0 = \Delta_f$; and
\item$\partial B_{i-i} \cap \partial B_i = \{x_i\}$ is a root of $B_i$ for $i \geq 1$.
\end{itemize}
The union
$$
\cR^B = \bigcup_{i=0}^\infty \overline{B_i}
$$
is called a {\it bubble ray}. The point $x_1 \in \partial \Delta_f$ is called the {\it root} of $\cR^B$. The accumulation set $\omega(\cR^B)$ of $\cR^B$ is defined as the accumulation set of the sequence $\{B_i\}_{i=0}^\infty$. Note that $\omega(\cR^B) \subset J_f$. If $\omega(\cR^B) = \{x_\infty\}$, then $x_\infty$ is called the {\it landing point} of $\cR^B$.

\begin{figure}[ht]
\includegraphics[width=1.1\textwidth, trim=100 0 0 0, clip]{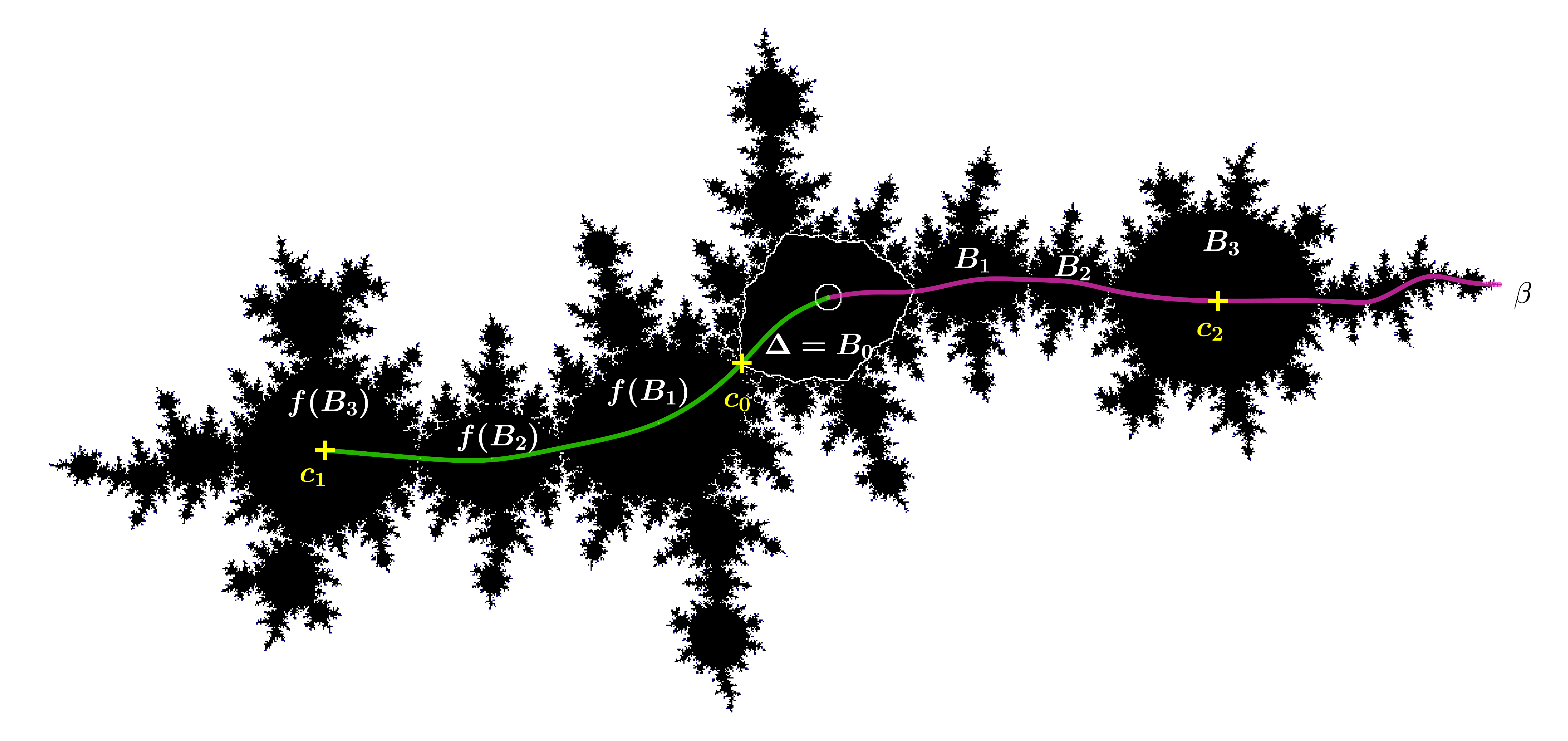}
\caption{The example of a quartic polynomial $f$ with the fixed Siegel disk $\Delta$ of bounded type rotation number. The critical points $c_1, c_2, c_3$ are marked with crosses. In this example, $c_0 \in \partial \Delta$, $f(c_2) = c_1$, and $f^3(c_1) \in \Delta$. The bubble ray $\mathcal R$ connecting $\Delta$ to the fixed point $\beta$ is marked with the purple curve. In this example, $f(\mathcal R) = \mathcal R \cup f(B_1) \cup f(B_2) \cup f(B_3)$, and hence $f(\mathcal R)$ passes through $\Delta$ and is not a bubble ray in our definition.}
\label{Fig:BubbleProblem}
\end{figure}

Note that in general, the image of a bubble ray is not necessarily a bubble ray in the sense of the definition given above (see Figure~\ref{Fig:BubbleProblem}). The image of a bubble ray contains a unique bubble ray. Nonetheless, the situation is more transparent for \emph{periodic bubble rays}. 

A bubble ray $\cR^B$ is {\it periodic} if $f^p(\cR^B) = \cR^B$ for some $p \geq 1$ and for every $k \in \{1, \ldots, p\}$, the image $f^k(\cR^B)$ is a bubble ray. We say that a bubble ray is {\it rational} if $f^n(\cR^B)$ is periodic for some $n \geq 0$ and for all $k \in \{1, \ldots, n-1\}$, the set $f^k(\cR)$ is a bubble ray. Note that all fixed bubble rays are rooted at $c_0$. A bubble ray $\cR^B$ is said to be {\it $d$-adic of generation $k\geq 0$} if $\cR^B$ is a rational bubble ray such that $f^k(\cR^B)$ is a fixed bubble ray, and $k$ is the smallest number for which this is true.

\begin{figure}
    \centering
    \includegraphics[trim=20 20 20 20, clip, width=1.0\textwidth]{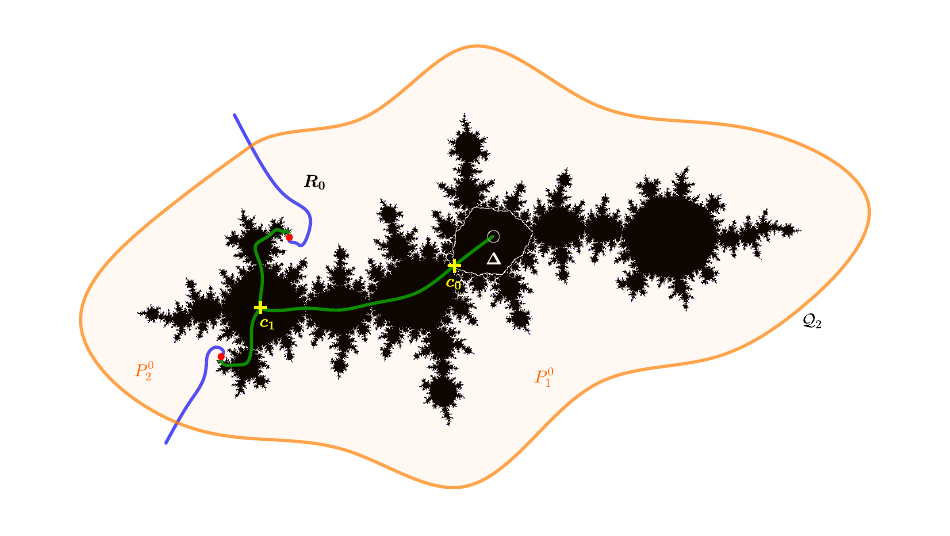}
    \caption{The construction of the initial puzzle for the same polynomial as in Figure~\ref{Fig:BubbleProblem}. The fixed bubble rays consist of the bubbles intersecting the green curves (the spines of the bubble rays). There are two fixed bubble rays in this example that overlap in the first four bubbles. The fixed external rays, shown in blue, land with the fixed bubble rays at the repelling fixed points (shown as red dots). The equipotential $\mathcal Q_2$ is indicated by the thick orange line. In this example, the initial puzzle partition consists of two puzzle pieces $P^0_1$ and $P^0_2$, shaded in orange.}
    \label{fig:Puzzle}
\end{figure}

A fixed bubble ray can be constructed in the following way. Starting with $B_0 = \Delta$, there exists at least one critical point $c_0 \in \partial \Delta$ (Theorem~\ref{siegel bound quasi}). Therefore, there exists at least one component $B_1 \in f^{-1}(B_0) \sm \{B_0\}$ that is attached to $B_0$ at $c_0$. Since $f(B_1) = B_0$, there exists at least one component $B_2$ in $f^{-1}(B_1) \sm \{B_1\}$ that is attached to $B_1$. Proceeding by induction, we construct a fixed bubble ray. A periodic bubble ray, say of period $p$, is constructed in the same way by passing to an iterate $f^p$ and constructing a fixed bubble ray for this iterate. Figure~\ref{fig:Puzzle} shows an example of two fixed bubble rays.

\begin{prop}[Rational bubble rays land]
\label{prop:land}
Every $p$-periodic bubble ray $\cR^B$ lands at a repelling or parabolic $p$-periodic point $x_\infty \in J_f$. \qed
\end{prop}

An {\it equipotential curve} at {\it level} $l \in (1, \infty)$ of $f$ is defined as
$$
\cQ_l := \{(\phi^\infty_f)^{-1}(le^{2\pi t i}) \; | \; t \in \bbR/\bbZ\}.
$$

Let $\bfR_0$ be the union of closures of all fixed bubble rays, all landing points of these bubble rays, and all external rays that also land at these points. Define the {\it initial puzzle} (or \emph{puzzle graph}) as
\begin{equation}\label{eq:0 puzzle}
\cZ_0 := \bfR_0 \cup \cQ_2.
\end{equation}

The {\it puzzle of depth} $n \geq 0$ is given by
$$
\cZ_n := f^{-n}(\cZ_0).
$$
A bounded connected component of $\bbC \setminus \cZ_n$ is called a {\it puzzle piece of depth $n$}. 

\begin{lem}[Puzzle pieces are Jordan domains]
    Every puzzle piece $P$ of depth at least $1$ is a Jordan domain.
\end{lem}

\begin{proof}
    The \emph{spine} of $f^{-n}(\bfR_0)$ is a finitely-branched tree $T_{-n} \subset f^{-n}(\bfR_0)$ whose root is at the Siegel fixed point and whose endpoints are all at infinity, such that for every external ray in $f^{-n}(\bfR_0)$, there exists a unique path from the root to the landing point of that ray that intersects each bubble along a Jordan arc.  

    Observe that for every $n \ge 1$, the tree $T_{-n}$ has at least two unbounded paths $\sigma_1, \sigma_2$ that are disjoint except at the root. In this way, $\ell := \overline{\sigma_1 \cup \sigma_2}$ contains the root and disconnects the plane. 
    
    It follows that each component of $\C \sm T_{-n}$ is a Jordan domain. Otherwise, there exists a component $V$ and a point $x \in \partial V$ that is doubly accessible from within $V$. Thus, there would exist open arcs $\gamma_1, \gamma_2 \subset V$ with common endpoints, one in $V$, and the other equal to $x$. In this way, $\overline{\gamma_1 \cup \gamma_2}$ bounds an open Jordan domain $V'$. Since the root of $T_{-n}$ is contained in $\ell$ and $\ell \cap V' = \emptyset$, the domain $V'$ does not contain the root. Consider an infinite path $\Gamma$ in $T_{-n}$ from the root that contains $x$ and intersects $V'$. Denote the subarc of $\Gamma$ from the root to $x$ by $\gamma$. Then $V' \supset \Gamma \sm \gamma$, which is a contradiction since $V'$ is bounded. 

    Now, let $P$ be a puzzle piece of depth $n \ge 1$. Each point $x \in\partial P$ belongs either to the spine, to the boundary of exactly one bubble $B$, or to the equipotential curve. In the first case, $x$ is uniquely accessible from within $P$ by the discussion above. In the last two cases, since $P$ does not contain $B$ or the points outside the equipotential curve, $x$ is again uniquely accessible. Since $\partial P$ is a piecewise Jordan curve (Theorem~\ref{siegel bound quasi}), it follows that $P$ is a Jordan domain.  
\end{proof}

The puzzle pieces are graded by their depth, and there are only finitely many puzzle pieces of any given depth. Furthermore, the puzzle pieces obey \emph{the Markov property}: 
\begin{itemize}
\item
every two puzzle pieces $P^n, P^m$, $n > m$, of depths $n$ and $m$ are either nested or disjoint; in the former case, $P^n \subset P^m$;
\item
$f(P^n) = P^{n-1}$ is a puzzle piece of depth $n-1$, and the map $f \colon P^n \to P^{n-1}$ is a proper holomorphic branched covering between two open topological disks.
\end{itemize}

The following theorem is a key technical result in the rigidity theory of Siegel polynomials.

\begin{thm}[Rigidity at Siegel boundary]
\label{Thm:RigidityBoundary}
Let $\{P^n\}_{n=0}^\infty$ be a sequence of puzzle pieces nested as $P^{n+1} \subset P^n$, and such that $\partial P^n \cap \partial \Delta_f \neq \emptyset$ for all $n \in \mathbb N$. Then
$$
\bigcap_{n=0}^\infty \overline{P^n} = \{z\} 
$$
for some $z \in \partial \Delta_f$. Consequently, the Julia set $J_f$ is locally connected at every point in $\partial \Delta_f$. \qed
\end{thm}

Combining Proposition~\ref{prop:land} and Theorem~\ref{Thm:RigidityBoundary}, we obtain:

\begin{thm}[Rigidity at Julia boundary]
\label{rigid boundary}
If $P$ is a puzzle piece for $f$, then $f$ is dynamically rigid at every point $z \in \partial P \cap J_f$ in the following sense: 

Let $Q^n$ be the union of closures of all puzzle pieces of depth $n \in \bbN$ containing the point $z$. Then 
\[
\bigcap \limits_{n=0}^\infty Q^n = \{z\}.
\] 
\qed
\end{thm}

It will be useful for us to use the terminology of fibers (see \cite{DS} for a discussion about this concept). Given a point $z \in J_f$, we define the \emph{fiber} $\fib(z)$ of $z$ to be the intersection of the closures of all puzzle pieces that contain $z$. The fiber of $z$ is \emph{trivial} if $\fib(z) = \{z\}$. In this way, Theorem~\ref{rigid boundary} says that all points at the boundary of each puzzle piece have trivial fibers.

\section{Background on complex box mappings}\label{sec:box}

In this section, following \cite{CDKvS}, we provide some background on the generalized renormalization technique via \emph{complex box mappings}.

\begin{defn}[Complex box mapping] 
\label{Def:BM}
A holomorphic map $F \colon \U \to \V$ between two open sets $\U \subset \V \subset \C$ is a \emph{complex box mapping} if the following holds:
\begin{enumerate}
\item
\label{It:DefBM1}
$F$ has finitely many critical points;
\item
\label{It:DefBM2}
$\V$ is the union of finitely many open Jordan disks with disjoint closures;
\item
\label{It:DefBM3}
every component $V$ of $\V$ is either a component of $\U$, or $V \cap \U$ is a union of Jordan disks with pairwise disjoint closures, each of which is compactly contained in $V$;
\item
\label{It:DefBM4}
for every component $U$ of $\U$ the image $F(U)$ is a component of $\V$, and the restriction $F \colon U \to F(U)$ is a proper map.
\item
\label{It:Nat1}
for each component $U$ of $\U$ there exists $n \ge 0$ so that $F^n(U) \setminus \U \neq \emptyset$;
\end{enumerate}
\end{defn}

\begin{rem}
In this paper, we slightly modify the definition of a complex box mapping given in, for instance, \cite{CDKvS}, by adding requirement \eqref{It:Nat1}. This condition guarantees that a box mapping has no cycles of components of $\U$ and is clearly satisfied if each component of $\U$ has escaping points. In \cite{CDKvS}, this condition was part of the dynamical naturality assumption and was called the \emph{no-permutation-condition}. We relax our definition of dynamically natural box mappings (see Definition~\ref{Def:NaturalBoxMapping}) by moving the aforementioned condition into the definition of a box mapping. 
\end{rem}

The \textit{filled Julia set of $F$} is defined as 
\[
K_F := \bigcap_{n \ge 0} F^{-n}(\V).
\]
The set $K_F$ consists of non-escaping points under the iteration of $F$, i.e., of points with orbits under $F$ that are defined for all iterates. 

For a given $n \ge 0$, a connected component of $F^{-n}(\V)$ is called a \textit{puzzle piece of depth $n$} for $F$. For $x \in \V$, the \textit{puzzle piece of depth $n$ containing $x$} is the connected component of $F^{-n}(\V)$ containing $x$. Note that if $x$ escapes $\U$ in fewer than $n$ steps, this set will be empty. 

Given a point $x \in \U$, the \emph{fiber of $x$} is the nested intersection of all puzzle pieces of $F$ containing $x$. We denote the fiber of $x$ as $\fib(x)$. By our convention, if the orbit of $x$ escapes, then $\fib(x) = \emptyset$.

We will denote the set of critical points for $F$ as $\Crit(F)$. Let us denote by $\PC(F)$ the union of the forward orbits of $\Crit(F)$. Note that $\PC(F)$ is different from the classical \emph{postcritical set}, which is defined as the closure of the union of the forward orbits of the critical values of $F$.  

\subsection{Dynamically natural complex box mappings}

Towards defining a more natural sub-family of complex box mappings, let us introduce the following two forward invariant subsets of $K_F$.

\begin{enumerate}
\item
Let $A \subset K_F$ be a finite set and $W$ be a union of finitely many puzzle pieces. We say that $W$ is a \textit{puzzle neighborhood} of $A$ if $A \subset W$ and each component of $W$ intersects the set $A$. If $\Crit(F) \neq \emptyset$, define 
\begin{equation*}
\begin{aligned}
\Koc(F) := \{x \in K_F \colon &\exists W \text{ puzzle neighborhood} \\ 
&\text{ of}\,\Crit(F) \text{ such that }\orb(x) \cap W = \emptyset\}
\end{aligned}
\end{equation*}
otherwise, i.e.\ when $F$ has no critical points, we set $\Koc(F) \equiv K_F$. Here, $\orb(x)$ stands for the orbit of $x$ under $F$.

\item
Let 
\[
m_F(x) := \modulus \left(\Comp_x \V \setminus \overline{\Comp_x \U}\right),
\] 
where $\modulus(\cdot)$ denotes the conformal modulus of an annulus, and $\Comp_x A$ is the connected component of a set $A$ containing $x$. For a given $\delta > 0$, set 
\[
K_\delta (F) := \left\{y \in K_F : \limsup_{k \ge 0} m_F(F^k(y)) > \delta \right\}.
\]
Define 
\[ 
\Kwi(F) := \bigcup_{\delta>0} K_\delta(F).
\] 
\end{enumerate}

\begin{defn}[Dynamically natural complex box mapping]
\label{Def:NaturalBoxMapping}
A complex box mapping $F \colon \U \to \V$ is called \emph{dynamically natural} if it satisfies the following conditions:
\begin{enumerate}
\item
\label{It:Nat2}
the Lebesgue measure of the set $\Koc(F)$ is zero;
\item
\label{It:Nat3}
$K_F = \Kwi(F)$.
\end{enumerate}
\end{defn}  

We refer the reader to \cite[Section 4.4]{CDKvS} for the discussion on the motivation for Definition~\ref{Def:NaturalBoxMapping}.

\subsection{Combinatorial equivalence of complex box mappings}

By Definition~\ref{Def:BM}, each component in the sets $\U$ and $\V$ for a box mapping $F\colon \U \to \V$ is a Jordan disk. Therefore, by the Carath\'eodory theorem, for every branch $F \colon U \to V$ of this mapping, there exists a well-defined continuous extension $\hat F \colon \cl U \to \cl V$, and by continuity, this extension is unique. (Here, $\cl A$ stands for the closure of a set $A$.) Let us denote by $\hat F \colon \cl \U \to \cl \V$ the total extended map. 

\begin{defn}[Itinerary of puzzle pieces relative to curve family]
\label{DefA:CurveFamily}
Let $F \colon \U \to \V$ be a complex box mapping, and let $X\subset \partial \V$ be a finite set with one point on each component of $\partial \V$. Let $\Gamma$ be a collection of simple curves in $(\cl \V) \setminus (\U \cup \PC(F))$, one for each $y \in \hat F^{-1}(X)$, that connects $y$ to a point in $X$. Then for every $n \ge 0$ and for each component $U$ of $F^{-n}(\U)$ there exists a simple curve connecting $X$ to $\partial U$ of the form $\gamma_0\ldots\gamma_n$ where $\hat F^k(\gamma_k) \in \Gamma$. The word $(\gamma_0, \hat F(\gamma_1), \ldots, \hat F^n(\gamma_n))$ is called the \emph{$\Gamma$-itinerary} of $U$.
\end{defn}

If a component $V$ of $\V$ is also a component
of $\U$, then the corresponding curve will be contained in the Jordan curve $(\cl V) \setminus (\U \cup \PC(F))
= \partial V$, because by definition $\PC(F)\subset \V$ and so no point in the boundary of $V$ lies
in $\PC(F)$. 

Note that the $\Gamma$-itinerary of $U$ is not uniquely defined even though there is a unique finite word for every $y'\in \hat F^{-n}(x)\cap \partial U$ (see \cite[Section 5.4]{CDKvS} for details). 

\begin{defn}[Combinatorial equivalence of box mappings]
\label{DefA:CombEquivBoxMappings}
Let $F \colon \U \to \V$ and $\tilde F \colon \tilde \U \to \tilde \V$ be two complex box mappings. Let $H \colon \V \to \tilde \V$ be a homeomorphism with the property that $H(\U)=\tilde \U$, $H(\PC(F)\setminus\U) = \PC(\tilde F) \setminus \tilde \U$ and such that it has a continuous extension $\hat H \colon \cl \V \to \cl \tilde \V$ to the closures of $\V$ and $\tilde \V$.

The maps $F$ and $\tilde F$ are called \emph{combinatorially equivalent w.r.t.\ $H$} if:
\begin{itemize}
\item
$H$ is a bijection between $\Crit(F)$ and $\Crit(\tilde F)$; for $c \in \Crit(F)$, $\tilde c := H(c)$ is the \emph{corresponding} critical point;
\item
there exists a curve family $\Gamma$ as in Definition~\ref{DefA:CurveFamily} so that for every $n \ge 0$, and for each $k \ge 0$ such that both $F^k(c)$ and $\tilde F^k(\tilde c)$ are well-defined, the $\Gamma$-itineraries of $\Comp_{F^k(c)} F^{-n}(\V)$ coincide with the $\hat H(\Gamma)$-itineraries of $\Comp_{\tilde F^k(\tilde c)} \tilde F^{-n}(\tilde \V)$.
\end{itemize} 
\end{defn}

\subsection{Renormalization of complex box mappings}

A complex box mapping $F$ is called \textit{(box) renormalizable} if there exists $s \ge  1$ and a puzzle piece $W$ of some depth containing a critical point $c \in \Crit(F)$ such that $F^{ks}(c) \in W$ for all $k \ge 0$. If this is not the case, then $F$ is called \emph{non-renormalizable}.

It is easy to see that for a complex box mapping, $F$ is box renormalizable if and only if it is polynomial-like renormalizable (see \cite[Lemma 5.1]{CDKvS}).

\subsection{Rigidity of non-renormalizable dynamically natural box mappings}

The following three theorems are the main results on rigidity of non-renormalizable complex box mappings. For the proofs, see \cite{KvS, CDKvS}. 

\begin{thm}[Dynamical rigidity]
\label{Thm:BoxDynRigidity}
If $F\colon \U\to \V$ is a non-renormalizable dynamically natural complex box mapping, then each point $x \in K_F$ is contained in arbitrarily small puzzle pieces. \qed
\end{thm}

Recall that for a forward-invariant subset $E \subset K(F)$, the \emph{line field} supported on $E$ is an assignment of straight lines $\mu(z)$ through the origin in the tangent space of each point $z \in E$ such that the slopes of $\mu(z)$ depend measurably on $z$. This line field is \emph{invariant under $F$} if $F'$ sends the line $\mu(z)$ to the line $\mu(F(z))$. In particular, if $K(F)$ has zero Lebesgue measure, then $K(F)$ does not support invariant line fields by definition. 

\begin{thm}[Absence of invariant line fields]
\label{Thm:BoxLineFields}
Let $F\colon \U\to \V$ be a non-renormalizable dynamically natural complex box mapping. Assume that each puzzle piece $P$ of $F$ contains an open set of escaping points, i.e.\ $\inter \left(P \,\sm K_F \right) \neq \emptyset$. Then $F$ carries no measurable invariant line fields supported on $K_F$. \qed 
\end{thm}

This theorem is a weaker version of \cite[Theorem~6.1(2)]{CDKvS}, where box mappings carrying invariant line fields are characterized, but it easily follows. This version will be sufficient for our purposes. 

\begin{thm}[Parameter rigidity]
\label{Thm:BoxParamRigidity}
Let $F\colon \U\to \V$ be a non-renormalizable dynamically natural complex box mapping. Suppose $\tilde F \colon \tilde \U \to \tilde \V$ is another dynamically natural complex box mapping for which there exists a quasiconformal homeomorphism $H \colon \V \to \tilde\V$ so that
\begin{enumerate}
\item
\label{part3c}
$\tilde F$ is combinatorially equivalent to $F$ w.r.t.\ $H$, and so in particular $H(\U) = \tilde \U$,
\item
\label{part3b}
$\tilde F \circ H = H \circ F$ on $\partial \U$, i.e.\ $H$ is equivariant on $\partial \U$.
\end{enumerate}
Then $F$ and $\tilde F$ are quasiconformally conjugate, and this conjugacy agrees with $H$ on $\V\setminus \U$.     \qed
\end{thm}

\section{Extracting dynamically natural complex box mapping}
\label{sec:DNBM}

Let $f$ be an atomic Siegel polynomial of bounded type with Siegel disk $\Delta_f$, and $\Crit^\star(f)$ be the set of critical points of $f$ that do not land\footnote{Here and everywhere, when we say ``a point lands in a set'' we mean ``the orbit of a point intersects a set''.} in $\overline \Delta_f$. In this section, we will associate to $f$ a dynamically natural complex box mapping that captures the dynamics of the points in $\Crit^\star(f)$. We will work with puzzle pieces constructed in Section~\ref{sec:puzzle}.

We start with some preliminary results and definitions. 

A set $A \subset \C$ is \emph{nice} if for every point $x \in \partial A$, one has that $f^n(x) \not \in \inter A$ for all $n \ge 1$. (Here, $\inter A$ stands for the interior of the set $A$.)

Given a holomorphic map $f \colon U \to V$ and a nice open set $A \subset U$, define 
\[
\DomE(A) := \{z \in U \colon \exists k \ge 1, f^k(z) \in A\}, \quad \Dom(A) := \DomE(A) \cap A, \quad \DomL(A):=\DomE(A) \cup A.
\]
The components of $\DomE(A)$, $\Dom(A)$ and $\DomL(A)$ are called called, respectively, first \emph{entry, return and landing domains (to $A$ under $f$)}. Note that $\cR(A) \subset \cL(A)$ and $\cR(A) = \cL(A) \cap A$. 

\begin{lem}[When a union of puzzle pieces is nice]
\label{Lem:UnionIsNice}
Let $\Y$ be a nice union of puzzle pieces of depths at most $k \ge 0$. If $Y$ is a puzzle piece of depth at least $k$ that is either a component of the first entry domain, or else is not contained in any of the components of the first entry domain to $\Y$, then $\Y \cup Y$ is a nice set.  

In particular, the union of disjoint puzzle pieces of the same depth is a nice set.
\end{lem}

(Note that in this lemma, $Y$ can be a component of the first entry domain.)

\begin{proof}[Proof of Lemma~\ref{Lem:UnionIsNice}]
Pick $x \in \partial(\Y \cup Y)$. To check that the union in question is nice, we distinguish two cases: 1) $x \in \partial \Y$ and 2) $x \in \partial Y$. 

In case 1), since $\Y$ is nice and open, $f^n(x) \not \in \Y$ for all $n \ge 1$. Furthermore, $f^n(x) \not \in Y$ for all $n \ge 1$ because otherwise, by the Markov property, a component of $\Y$ that has $x$ on its boundary would be of depth at least $k+1$, which is impossible by assumption.  

In case 2), $f^n(x) \not \in Y$ for all $n \ge 0$ because $Y$ is nice and open, and $f^n(x) \not \in \Y$ for all $n \ge 0$ because $Y$ does not lie in any of the landing domains to $\Y$ and hence $Y$ is never mapped into a component of $\Y$.
\end{proof}

The following statements are useful in constructing complex box mappings.

\begin{lem}[First return to protected puzzle pieces]
\label{Lem:FRM}
Let $\X, \Y$ be two nice unions of finitely many puzzle pieces such that $\X \Subset \Y$. Suppose that each component $X$ of $\X$ has the following \emph{disjointness property}:
\begin{itemize}
\item
if there exists $k \ge 1$ such that $\partial f^k(X) \cap \Y \neq \emptyset$, the puzzle pieces $f(X), \ldots, f^{k}(X)$ are disjoint from $\X$.
\end{itemize}
If $\X'$ is a union of components of $\cL(\X) \sm \X$, then each component of $\cR(\X \cup \X')$ is either a component of $\X'$, or is compactly contained in $\X$.
\end{lem}

\begin{proof}
First of all, observe that the union $\X \cup \X'$ is a nice set (by Lemma~\ref{Lem:UnionIsNice}). Hence, the domain of the first return map to $\X \cup \X'$ is a finite union of puzzle pieces. Further note that $\cR(\X \cup \X') \cap \X = \cR(\X)$, and hence it is enough to show that each component of $\cR(\X)$ is compactly contained in $\X$.

Let $P$ be a component of $\cR(\X)$, and $X$ be the component of $\X$ that contains $P$. Our goal is to show that $P \Subset X$. 

Assume the contrary, i.e., $\partial P \cap \partial X \neq \emptyset$. If $k \ge 1$ is the first return time of $P$ to $\X$, then $f^k(P)=:X'$ is a component of $\X$ and $X' \subset f^k(X)$. Since $\partial P \cap \partial X\neq \emptyset$, the same is true for $X'$ and $f^k(X)$, i.e., $\partial X' \cap \partial f^k(X) \neq \emptyset$. By assumption, $X' \Subset \Y$, and hence $f^k(X) \subset \Y$ (by the Markov property) and the boundary of $f^k(X)$ intersects $\Y$. Because $f^k(X')$ contains $X'$, the set $f^k(X)$ intersects $\X$, and this is a contradiction to the disjointness property. 
\end{proof}

Now we need the following lemmas about the weak protection of critical puzzle pieces. Recall that a critical point $c \in \Crit^\star(f)$ is called \emph{combinatorially recurrent} if its orbit intersects every puzzle piece around $c$. Note that, since the rational lamination of $f$ is empty, $\Crit^\star(f)$ contains no renormalizable critical points, and hence critical points can be recurrent only in a non-periodic way. Let us distinguish the following subsets of $\Crit^\star(f)$: the first one, which we call $\Crit^\star_r(f)$, will be a set of all combinatorially recurrent critical points, while $\Crit^\star_{nr}(f)$ will stand for all combinatorially non-recurrent critical points. Obviously,
\[
\Crit^\star(f) = \Crit_r^\star(f) \sqcup \Crit_{nr}^\star(f). 
\]  

We say that a puzzle piece $P$ is \emph{weakly protected} if there exists a puzzle piece $Q$ such that $P \Subset Q$. 

Let $P_n \supset P_{n+k}$ be two nested puzzles of depth $n$ and $n+k$ respectively. The puzzle piece $P_{n+k}$ is referred to as a {\it child of $P_n$} if
$$
f^i(P_{n+k}) \cap P_{n+k} = \varnothing
\matsp{for}
0 < i < k
\matsp{and}
f^k(P_{n+k}) = P_n.
$$

\begin{lem}[Two Children Lemma]\label{two kids}
For every $c \in \Crit_r^\star(f)$ and $n \geq 0$, there exist at least two values $0<  k_1 < k_2$ such that $P_{n+k_1}(c)$ and $P_{n+k_2}(c)$ are children of $P_n(c)$.
\end{lem}

\begin{proof}
The first child $P_{n+k_1}(c)$ is given by the component of the first return domain to $P_n(c)$ containing $c$. Since $f$ is not renormalizable, there exists a smallest number $m \in \bbN$ such that $f^{mk_1}(c) \notin P_{n+k_1}(c)$. Let $l \in \bbN$ be the smallest number such that $f^{mk_1+l}(c) \in P_n(c)$. We claim that the piece $P_{n+k_2}(c)$ with $k_2 := mk_1+l$ is another child of $P_n(c)$.

Suppose towards a contradiction that $f^i(P_{n+k_2}(c)) = P_{n+k_2-i}(c)$ for some $i < k_2$. If $i < l$, then mapping
$$
P_{n+k_2}(c) \subsetneq P_{n+k_2-i}(c)
$$
by $f^{mk_1}$, we obtain
$$
P_{n+l}(f^{mk_1}(c)) \subsetneq P_{n+l-i}(f^{mk_1}(c)).
$$
Note that $P_{n+l-i}(f^{mk_1}(c)) \subsetneq P_n(c)$ since $f^{(m-1)k_1}(c) \in P_{n+k_1}(c)$. On the other hand, $f^{l-i}$ maps $P_{n+l-i}(f^{mk_1}(c))$ to $P_n(c)$. This contradicts the fact that $P_{n+l}(f^{mk_1}(c))$ is the component of the first return domain of $P_n(c)$ containing $f^{mk_1}(c)$.

Suppose that
$$
P_{n+jk_1}(c) \subsetneq P_{n+k_2-i}(c) \subsetneq P_{n+(j-1)k_1}(c)
$$
for some $1 \leq j \leq m$. Mapping by $f^{(j-1)k_1}$, we obtain
$$
P_{n+k_1}(c) \subsetneq f^{(j-1)k_1}(P_{n+k_2-i}(c)) \subsetneq P_n(c).
$$
This implies that $f^{(j-1)k_1}(P_{n+k_2-i}(c)) = P_{n+k'}(c)$ with
$$
k' = (m-j+1)k_1+l-i < k_1,
$$
which must map onto $P_n(c)$ under $f^{k'}$. This contradicts the fact that $P_{n+k_1}(c)$ is the first child of $P_n(c)$. Thus, $i = k_2 - jk_1$ for some $1 \leq j \leq m$.

Note that
$$
f^i(P_{n+k_2}(c)) = P_{n+jk_1}(f^i(c)) = P_{n+jk_1}(c).
$$
Since $f^{(j-1)k_1}(c) \in P_{n+k_1}(c)$, it follows that
$$
f^{i+(j-1)k_1}(P_{n+k_2}(c)) = P_{n+k_1}(c).
$$
Observe that $i+(j-1)k_1 = k_2 - k_1 = (m-1)k_1 +l$. Mapping by $f^{k_1-l}$, we obtain
$$
f^{mk_1}(c) \in f^{mk_1}(P_{n+k_2}(c)) = f^{k_1-l}(P_{n+k_1}(c)) = P_{n+l}(c) \subsetneq P_n(c).
$$
Note that $P_{n+l}(c)$ maps to $P_n(c)$ under $f^l$ with $l < k_1$. This again contradicts the fact that $P_{n+k_1}(c)$ is the first child of $P_n(c)$.
\end{proof}

\begin{thm}[Recurrent critical points are weakly protected]
\label{crit protect}
For every $c \in \Crit_r^\star(f)$ and $N \in \bbN$, there exists $n \geq N$ and $k \in \bbN$ such that  $P_{n+k}(c) \Subset P_n(c)$, and $P_{n+k}(c)$ is a child of $P_n(c)$.
\end{thm}

\begin{proof}
We simplify notation by denoting $P_n(c) = P_n$. By \lemref{two kids}, there exist a nested sequence
$$
P_0 =: P_{m_0} \supsetneq P_{m_1} \supsetneq \ldots
$$
of successive second children (i.e. $P_{m_{i+1}}$ is the second child of $P_{m_i}$). We will show that
$$
m_{i+1}-m_i > m_i
\matsp{for}
i \in \bbN.
$$

Let $P_{n_{i+1}}$ be the first child of $P_{m_i}$. Observe that
$$
P_{m_{i+1}} \subsetneq P_{n_{i+1}} \subsetneq P_{m_i}.
$$
Applying $f^{a_{i+1}} := f^{n_{i+1}-m_i}$:
\begin{equation}
    \label{Eq:Chain1}
f^{a_{i+1}}(P_{m_{i+1}}) \subsetneq f^{a_{i+1}}(P_{n_{i+1}}) = P_{m_i} \subsetneq P_{m_{i-1}} \subseteq f^{a_{i+1}}(P_{m_i}).    
\end{equation}
Note that the last inclusion follows from the fact that, since $P_{m_{i}}$ is a child of $P_{m_{i-1}}$, the orbit of $P_{m_{i}}$ cannot contain a critical point (in particular, cannot contain $P_{m_i}$) before returning to $P_{m_{i-1}}$. Hence, the inclusion $P_{m_i} \subset f^{a_{i+1}}(P_{m_i}) \subset P_{m_{i-1}}$ is impossible.

The chain \eqref{Eq:Chain1} implies that $a_{i+1}$ is strictly less than the return time $m_{i+1}-m_i$ of $P_{m_{i+1}}$, while it is greater than the return time $m_i - m_{i-1}$ of $P_{m_i}$.

Subsequently, apply $f^{a_i}$:
$$
f^{a_{i+1}+a_i}(P_{m_{i+1}}) \subsetneq f^{a_i}(P_{m_i})\subsetneq f^{a_i}(P_{n_i}) = P_{m_{i-1}} \subsetneq P_{m_{i-2}} \subseteq f^{a_i}(P_{m_{i-1}}) \subseteq f^{a_{i+1}+a_i}(P_{m_i}).
$$
This means that $a_{i+1}+a_i$ is strictly less than $m_{i+1}-m_i$. Indeed, assume the contrary, namely, $a_{i+1}+a_i \ge m_{i+1}-m_i$; then $t :=  a_{i+1}+a_i - (m_{i+1}-m_i) \ge 0$ measures the time between the moment the orbit of $P_{m_{i+1}}$ lands at $P_{m_i}$ and before it landed at $f^{a_{i+1}+a_i}(P_{m_{i+1}}) = f^t(P_{m_i})$. Since $f^{a_{i+1}+a_i}(P_{m_{i+1}}) \subsetneq f^{a_i}(P_{m_i})$, we have $t < a_i$. On the other hand, by \eqref{Eq:Chain1}, $a_{i+1} < m_{i+1} - m_i$, hence, 
\[
a_i - t = a_{i+1}-(m_{i+1}-m_i) < 0,
\]
which is a contradiction.

Proceeding by induction, we conclude that
$$
A_{i+1} := \sum_{l=1}^{i+1} a_l
$$
is strictly less than $m_{i+1}-m_i$, while $P_0 \subseteq f^{A_{i+1}}(P_{m_i}).$ Thus, $A_{i+1} \geq m_i$, and the claim follows.

By \thmref{rigid boundary}, there exist $n_0 \geq 1$ such that $P_n \Subset P_0$ for all $n \geq n_0$. Choose $i \in \bbN$ sufficiently large so that $m_i \geq n_0$. By the claim, the puzzles $P_{m_{i+1}} \subset P_{m_i}$ map under $f^{m_{i+1}-m_i}$ to $P_{m_i}$ and $f^{m_{i+1}-m_i}(P_{m_i}) \supset P_0$. Consequently, by pulling back the annulus
$$
f^{m_{i+1}-m_i}(P_{m_{i+1}}) = P_{m_i} \Subset P_0
$$
under $f^{-m_{i+1}+m_i}$, we see that $P_{m_{i+1}} \Subset P_{m_i}$. This construction can be done at any sufficiently deep level $m_i \geq N$.
\end{proof}

We will now use this result to show how to extract a complex box mapping associated to an atomic Siegel polynomial. The main result is the following theorem.

\begin{thm}[Extracting Siegel box mapping]
\label{Thm:Extract}
Let $f$ be a non-renor\-malizable bound\-ed-type Siegel polynomial, and let $\Crit^\star(f) \subset \Crit(f) \cap J_f$ be the set of critical points of $f$ that do not land in the boundary of Siegel disks. Then there exists a complex box mapping $F \colon \U \to \V$ such that:
\begin{enumerate}
\item
the puzzle pieces of $F$ are puzzle pieces of $f$;
\item
$\Crit(F) \subset \Crit^\star(f) \subset \V$;
\item
if the orbit of a point $z \in J_f$ intersects $\V$, then this orbit intersects $K(F)$;
\item
$F$ satisfies the condition $\Kwi(F) = K_F$ from the definition of dynamical naturality (Definition~\ref{Def:NaturalBoxMapping}). 
\end{enumerate}
\end{thm}

\begin{proof}
The required box mapping will be constructed as a first return map to a carefully chosen nice puzzle neighborhood of $\Crit^\star(f)$. We construct this neighborhood by induction on the cardinality $\#\Crit_r^\star(f)$ of the set $\Crit_r^\star(f)$ of the recurrent critical points in $\Crit^\star(f)$. It will be convenient to start the whole construction at sufficiently large depth $s \ge 0$ chosen so that 
\begin{equation}
\label{Eq:Ass1}
\text{if  }\orb(c) \cap P_s(c') \neq \emptyset \text{ for some }c, c' \in \Crit^\star(f), \text{ then }\overline{\orb}(c) \cap \fib(c') \neq \emptyset.
\end{equation}  
In words, this means that if the orbit of a critical point $c$ enters a critical puzzle piece of depth $s$ around a critical point $c'$, then the orbit of $c$ must accumulate on the fiber of $c$. Such a choice is possible because we have only finitely many critical points. Note that it is possible that $c = c'$. In particular, if $c \in \Crit_{nr}^\star(f)$, then the orbit of $c$ never re-enters $P_s(c)$. Furthermore, by increasing $s$ if necessary, we can assume that 
\begin{equation}
\label{Eq:Ass2}
\text{if } c' \in P_s(c), \text{then }c' \in \fib(c).
\end{equation}
In particular, the last assumption implies that $\fib(c) = \fib(c')$.

In order to construct the required neighborhood of $\Crit^\star(f)$, we inductively define a sequence of sets $\cX_0 = \cY_0 = \emptyset, \cX_1 \Subset \cY_1, \ldots, \cX_m \Subset \cY_m$, $m \le \# \Crit^\star_r(f)$ such that each pair $\cX_i, \cY_i$ satisfies the assumption of Lemma~\ref{Lem:FRM} and the only components of $\cX_i$ and $\cY_i$ are critical puzzle pieces of depth at least $s$. Suppose $\cX_k \Subset \cY_k$ are defined. If all critical points in $\Crit_r^\star(f) \sm \cY_k$ eventually land into $\cY_k$, then we end the induction with the sets $\cX' := \cX_k$ and $\cY':=\cY_k$. Otherwise, choose a critical point $c_{k+1} \in \Crit_r^\star(f)$ whose orbit never intersects $\cX_k$. By \eqref{Eq:Ass1} and \eqref{Eq:Ass2}, this implies that the orbit of $c_{k+1}$ never intersects $\cY_{k+1}$. Let $Y_{k+1} \Supset X_{k+1} \ni c_{k+1}$ be a pair of puzzle pieces given by Theorem~\ref{crit protect}. By that theorem, we can choose $Y_{k+1}$ to be disjoint from $\cY_{k}$ and of depth larger than the depths of the components of $\cX_{k}$ (and hence of $\cY_{k}$). In this way, $\cY_{k} \cup Y_{k+1}$ and $\cX_{k} \cup X_{k+1}$ are nice set. Indeed, the orbit of $c_{k+1}$ never enters $\cY_{k}$ and $\cX_{k}$, and hence $Y_{k+1}$, respectively $X_{k+1}$ cannot be component of the first entry domain to $\cY_{k}$, respectively $\cX_{k}$; the niceness now follows by Lemma~\ref{Lem:UnionIsNice}.

Finally, define $\cY_{k+1} := \cY_k \cup Y_{k+1}$ and $\cX_{k+1} := \cX_{k} \cup X_{k+1}$. Let us check that $\cX_{k+1}$ and $\cY_{k+1}$ satisfy the assumptions of Lemma~\ref{Lem:FRM}. Clearly, $\cX_{k+1} \Subset \cY_{k+1}$. Let us check the disjointness property. By construction, since the orbit of $c_{k+1}$ never lands in $\Y_k$, the orbit of $X_{k+1}$ has a trivial disjointness property with respect to the components of $\cY_k$, namely, $\partial f^k(X_{k+1}) \cap \cY_k = \emptyset$ for all $k \in \bbN$. On the other hand, the orbit of $X_{k+1}$ has the disjointness property with respect to $Y_{k+1}$ by Theorem~\ref{crit protect}. Note that the orbit of the boundary of $\cX_k$ never intersects $Y_{k+1}$ by our choice of depths. Therefore, $\cX_{k+1}$, $\cY_{k+1}$ satisfy the assumptions of Lemma~\ref{Lem:FRM}. The induction step is justified, and the induction ends with sets $\cX' \Subset \cY'$.

By construction, each point in $\Crit^\star_r(f)$ is either contained in, or is mapped to $\cX'$. Some of the non-recurrent critical points $\Crit^\star_{nr}(f)$ can also map to $\cX'$. Now let $\Crit^{\circ}_{nr}(f)$ be the remaining non-recurrent critical points whose orbits do not intersect $\cX'$ (and thus, by assumption \eqref{Eq:Ass1}, also do not intersect $\cY'$). For each such critical point $c$, let us choose a puzzle piece $c \in Q_c$ of depth $\ell + 1$, where $\ell$ is the largest depth of puzzle pieces in $\cY'$ (note that $\ell > s$). Furthermore, since the points on the Julia boundary of Siegel puzzle pieces have trivial fibers (Theorem~\ref{rigid boundary}), for every $Q_c$ we can choose a puzzle piece $P_c$ such that $c \in P_c \Subset Q_c$, i.e., $P_c$ is weakly protected by $Q_c$. 

Let us now define
\[
\cX := \cX' \cup \bigcup_{c \in \Crit^{\circ}_{nr}(f)} P_c, \quad \cY := \cY' \cup \bigcup_{c \in \Crit^{\circ}_{nr}(f)} Q_c
\]
  
By construction, $\cX \Subset \cY$. By Lemma~\ref{Lem:UnionIsNice}, $\cY$ is a nice set. Finally, $\X$ satisfies the assumption of Lemma~\ref{Lem:FRM}, and hence if we define $\V$ to be $\cX$ union all the critical components of $\cL(\cX) \sm \cX$, then by Lemma~\ref{Lem:FRM}, the first return map $F \colon \U \to \V$ to $\V$ has the desired structure of a complex box mapping. The asserted properties (1)--(3) follow by construction. Assertion (4) follows verbatim as in Step 2 of the proof of \cite[Lemma 5.1]{DS} using the fact that each puzzle piece in $\cX$ is weakly protected by a puzzle piece in $\cY$.
\end{proof}

Let us now show that the box mapping constructed in Theorem~\ref{Thm:Extract} is dynamically natural. For that, we need to check the condition that the Lebesgue measure $\Koc(F)$ is zero (see Definition~\ref{Def:NaturalBoxMapping}). 

\begin{prop}
\label{Prop:NAT}
The box mapping constructed in Theorem~\ref{Thm:Extract} is dynamically natural in the sense of Definition~\ref{Def:NaturalBoxMapping}. 
\end{prop}

This proposition is an immediate corollary to the following much stronger statement:

\begin{thm}[Zero measure of accumulation on Siegel disk]
\label{Thm:AccZero}
Let $f$ be a Siegel polynomial. Then the set 
\[
J_\Delta := \{x \in J_f \colon \omega(x) \subset \partial \Delta_f\}
\]
has zero Lebesgue measure.
\end{thm}

\begin{proof}[Proof of Proposition~\ref{Prop:NAT}]
By construction of Siegel puzzle pieces, $\Koc(F) \subset J_\Delta$. The claim follows.
\end{proof}

\begin{proof}[Proof of Theorem~\ref{Thm:AccZero}]

The proof is based on the Lebesgue density point theorem and is similar to \cite[Lemma~3.12]{W23}. We repeat the proof for the sake of completeness. The main ingredient is the following key theorem from \cite{McM98, W23}:

\begin{thm}
\label{Thm:WMcM}
For every point $x \in J_\Delta$ and scale $r > 0$, there is a univalent map between pointed disks 
\[
f^i \colon (U, y) \to (V, c)
\]
of bounded distortion, where $c \in \partial\Delta_f$ is a critical point, the distance between $y$ and $x$ is of order $r$, and $(U,y)$ is a pointed disk of bounded shape with diameter of order $r$.  \qed
\end{thm}

(Recall that a pointed disk $(U,y)$ has bounded shape if there exists a constant $C > 0$ (that depends only on $f$) such that $\partial U$ contains the disk of radius $t$ centered at $y$ and is contained in the disk of radius $T$, also centered at $y$, with $T/t \le C$.)

Pick a point $x \in J_\Delta$ and sufficiently small $r > 0$. We will show that the round disk $\mathbb D(x, r)$ contains a smaller disk $D$ of radius comparable to $r$ such that $D \subset F_f$. The claim follows from the Lebesgue Density Point Theorem.

By Theorem~\ref{Thm:WMcM}, there is a univalent map $f^i \colon (U, y) \to (V, c)$ and the disk $(U,y)$ is contained in $\mathbb D(x,r)$. Since $\Delta_f$ is a quasi-disk, there is a round disk $D'$ in $V \cap \Delta_f \subset F_f$ of diameter comparable to that of $V$. Using the bounded distortion property, we can now lift the disk $D'$ to $U$. This will produce a disk (maybe smaller) of radius comparable to $r$ that lies in $\mathbb D(x, r) \cap F_f$.
\end{proof}

\section{Proof of Main Theorem~\eqref{it:main:1}, \eqref{it:main:2}}
\label{sec:proofAB}

Recall that if $f$ is an atomic Siegel polynomial, then the Julia set of $f$ is connected. Furthermore, $f$ has a unique Siegel disk $\Delta_f$ centered at $0$ and has no other neutral or attracting periodic points in $\C$. 

In this section, we prove Main Theorem~\eqref{it:main:1}, \eqref{it:main:2} in the following form:

\begin{thm}
\label{Thm:DynRigidity}
Let $f$ be an atomic Siegel polynomial of bounded type. Then
\begin{enumerate}
\item
all fibers of points in the Julia set $J_f$ are trivial, and hence, $J_f$ is locally connected;
\item
$J_f$ carries no invariant line fields.
\end{enumerate}
\end{thm}

For $f$ as in the theorem above, let $F \colon \U \to \V$ be the box mapping constructed in Theorem~\ref{Thm:Extract}. We have $\Crit^\star(f) \subset \V$. Recall that $\Crit^\star(f)$ is the set of critical points of $f$ in $J_f$ that do not lie in the puzzle graph of any depth. Let $\fs$ be the smallest depth of the puzzle graph that contains all points in $\left(\Crit(f) \cap J_f\right) \sm \Crit^\star(f)$. If the latter set is empty, we set $\fs=0$.

\begin{proof}[Proof of Theorem~\ref{Thm:DynRigidity}]
We start by proving the first assertion. 

Let $z \in J_f$ be a point. If the orbit of $z$ eventually lands in $\partial \Delta_f$, then $\fib(z) = \{z\}$ by \cite{Y}. So we can assume that $z$ does not lie on the puzzle graph of any depth. 

If the orbit of $z$ accumulates on a point in $\Crit^\star(f)$, and hence intersects $\V$ infinitely many times, then by assertion $(3)$ of Theorem~\ref{Thm:Extract}, this orbit intersects $K(F)$. Each point in $K(F)$ has trivial fiber, hence $\fib(z) = \{z\}$. 

The only case left to consider is if the orbit of $z$ does not accumulate on $\Crit^\star(f)$, and hence there exists a puzzle neighborhood $X$ of $\Crit^\star(f)$ so that $\orb(z) \cap X = \emptyset$.  We can assume that each puzzle piece in $X$ is of depth at least $\fs$.

There are three cases:
\begin{enumerate}
\item[a)]
there exists a point $w \in \omega(z)$ that does not belong to the puzzle graph of any depth.
\item[b)]
$\omega(z)$ lies in the puzzle graph and intersects the repelling fixed points that are the landing points of the Siegel bubble rays used in the construction of the initial puzzle partition.
\item[c)]
$\omega(z)$ only intersects the grand orbit of $\partial \Delta_f$.
\end{enumerate}

In case a), the idea is standard: there is a pair $P_m$, $P_{m+\ell}$ of puzzle pieces of depths $m > \fs$ and $m+\ell$ so that $w \in P_{m+\ell} \Subset P_m$. Since the orbit of $z$ accumulates on $w$, we can find a sequence $(k_i)$ so that $f^{k_i}(z) \in P_{m+\ell}$. For each $k_i$, we can pull back conformally the annulus $A:= P_{m} \sm \overline{P}_{m+\ell}$ by $f^{k_i}$ to obtain a sequence of annuli $A_i := P_{m+k_i}(z) \sm \overline{P}_{m+\ell+k_i}(z)$ around $z$. By conformality, $\mod(A_i) = \mod(A)$. Since these annuli can be assumed to be nested, the Gr\"otzsch inequality yields $\fib(z) = \{z\}$. 

In case b), the claim follows as in \cite[Lemma 4.7]{DS} using the fact that the accumulation points are repelling.

In case c), we can use the same argument as in a): Since the orbit of $z$ accumulates only on $\partial \Delta_f$, it also accumulates on $c \in \partial \Delta_f$. By Theorem~\ref{Thm:RigidityBoundary}, we can find a sequence of nested puzzle \emph{disks} $Q_1 \Supset Q_2 \Supset \ldots$ containing $c$ and shrinking to it \cite{Y}. We can now repeat the argument from above. 

Finally, for assertion (2), we decompose $J_f$ into 3 disjoint sets $A_1, A_2$, and $A_3$ as follows: $z \in A_1$ if $\omega(z) \cap K(F) \neq \emptyset$; $z \in A_2$ if $\orb(z)$ does not intersect some puzzle neighborhood of $\Crit^\star(f)$ but $\omega(z) \subset \partial \Delta_f$; and $A_3 := J_f \sm (A_1 \cup A_2)$. 
\begin{itemize}
\item
The set $A_3$ has zero Lebesgue measure by a standard argument (the set of points in $J_f$ that do not accumulate on the recurrent critical points has measure zero). Hence, $A_3$ supports no invariant line fields.
\item
The set $A_2$ has zero Lebesgue measure by Theorem~\ref{Thm:AccZero}, and hence again supports no invariant line fields.
\item
And finally, if $z \in A_1$, then $\orb(z) \cap K(F) \neq \emptyset$ (assertion $(3)$ of Theorem~\ref{Thm:Extract}). By Proposition~\ref{Prop:NAT} and Theorem~\ref{Thm:BoxLineFields}, $K(F)$ does not support invariant line fields. Hence, so does $A_1$. 
\end{itemize}

Altogether, we obtain that $J_f$ supports no invariant line fields.
\end{proof}

\section{Oriented Combinatorial Trees}
\label{sec:combtrees}

A (combinatorial) tree $T = (\cV, \cE, v^{\rt})$ is given by its vertex set $\cV$, its edge set $\cE$ and its root $v^{\rt} \in \cV$. We allow $T$ to contain self-loops: that is, there may be an edge $e \in \cE$ that connects a vertex $v \in \cV$ to itself (i.e. $e = (v,v)$). For $v \in \cV$, denote by $\cE_v \subset \cE$ the set of edges that connect to $v$. We say that $T$ is {\it oriented} if every vertex carries an orientation: a fixed cyclic order on $\cE_v$ for all $v \in \cV$. For any subset $X \subset T$, define
$$
\cE_X := \bigcup_{v \in X \cap \cV} \cE_v.
$$

A path in $T$ is a set given by
$$
\Gamma = v_0 \cup \bigcup_{i=1}^K (e_i \cup v_i)
\matsp{for some}
K \in \bbN \cup \{0, \infty\}
$$
such that $e_i = (v_{i-1}, v_i)$. The length of $\Gamma$ is denoted by $|\Gamma| = K$.

We say that a path $\Gamma$ is {\it rooted} if it is minimal (i.e. no loops) and starts at $v_0 = v^{\rt}$. For $v \in V$, let $\Gamma_T(v)$ be the unique rooted path containing $v$. The {\it height of $v$} is defined as
$$
\het(v) := |\Gamma_T(v)|-1.
$$
If $v \neq v^{\rt}$, the {\it joint of $v$} is the unique edge $\jt(v) \in \cE_v$ contained in $\Gamma_T(v)$. In this case, we define $\het(\jt(v)) := \het(v)$.

Let $\Gamma$ be a rooted path of length $0 \leq K \leq \infty$. Then there are ordered sequences $\{v^i\}_{i=0}^K$ and $\{e^i\}_{i=0}^{K-1}$ of vertices and edges whose union is equal to $\Gamma$ such that $v^0 = v^{\rt}$ and $e^i = (v^i, v^{i+1}) = \jt(v^{i+1})$ for $0 \leq i < K$. We write $\Gamma \sim \{v^i\}_{i=0}^K.$ The set $\{e^i\}_{i=0}^{K-1}$ is referred to as the {\it joint set of $\Gamma$}.

\begin{defn}
Let $\tiT = (\ticV, \ticE, \tiv^{\rt})$ and $T=(\cV, \cE, v^{\rt})$ be oriented trees. A {\it tree homomorphism} $g : \tiT \to T$ is a map such that
\begin{enumerate}
\item for $v \in \ticV$, we have $g(v) \in \cV$;
\item for $e = (v,u) \in \ticE$, we have $g(e) = (g(v), g(u)) \in \cE$;
\item $g$ is orientation-preserving on each vertex; and
\item $g(\tiv^{\rt}) = v^{\rt}$.
\end{enumerate}
If $g$ is bijective, then it is called a {\it tree isomorphism}.
\end{defn}

\begin{defn}
A {\it  tree cover $g : \tiT \to T$} is a tree homomorphism such that the following properties are satisfied.
\begin{enumerate}
\item For $v \in \ticV$, there exists $\deg(v) \geq 1$, referred to as the {\it degree of $v$}, such that every $e \in \cE(g(v))$ has exactly $\deg(v)$ preimages in $\ticE_v$.
\item The set $\ticV^{\crit}(g) \subset \ticV$ of vertices whose degrees are at least $2$ is finite.
\item For $e = (v, u) \in \ticE$, we have $g(v) \neq g(u)$.
\end{enumerate}
A vertex $c \in \ticV^{\crit}(g)$ is said to be {\it critical}. Define $\tiT^{\crit}(g)$ as the minimal subtree of $\tiT$ containing $\ticV^{\crit}(g)$ (by convention, we set $T^{\crit}(g) := \{\tiv^{\rt}\}$ if $\ticV^{\crit}(g) = \varnothing$).
\end{defn}

\begin{rem}
A  tree cover without a critical vertex is a tree isomorphism.
\end{rem}

\begin{prop}\label{cov ext}
Let $g_1 : \tiT_1 \to T$ and $g_2 : \tiT_2 \to T$ be  tree covers. If there exists a bijection
$$
\phi_{\crit} : \tiT^{\crit}(g_1) \cup \ticE_{\tiT^{\crit}(g_1)} \to \tiT^{\crit}(g_2) \cup \ticE_{\tiT^{\crit}(g_2)}
$$
that restricts to a  tree isomorphism from $\tiT^{\crit}(g_1)$ to $\tiT^{\crit}(g_2)$ such that $g_1 = g_2 \circ \phi_{\crit}$, then $\phi_{\crit}$ extends uniquely to a  tree isomorphism $\phi : \tiT_1 \to \tiT_2$ such that $g_1 = g_2 \circ \phi.$
\end{prop}

\begin{proof}
If $\tiV^{\crit}_1(g_1)=\tiV^{\crit}_2(g_2)=\varnothing$, then $g_1$ and $g_2$ are  tree isomorphisms, and $\phi := g_2^{-1} \circ g_1$. Uniqueness follows from the fact that the only tree isomorphism $\phi : T \to T$ that fixes $\cE_{v^{\rt}}$ is the identity (recall that $\phi$ must preserve orientation at every vertex).

Proceeding by induction, suppose that the result holds for tree covers with at most $n$ critical vertices. Suppose that $g_1$ and $g_2$ have $n+1$ critical vertices. We modify these covers to reduce the number of critical vertices.

Choose a critical vertex $c_1 \in \tiT^{\crit}(g_1)$ with maximal height. Let $d := \deg(c_1) \geq 2$. Denote the joint at $c_1$ by $\jt(c_1)$. The corresponding critical vertex for $g_2$ is given by $c_2 := \phi_{\crit}(c_1)$ with joint $\jt(c_2)$. Denote $\nu := g_1(c_1)= g_2(c_2)$ and $\jt(\nu) := g_1(\jt(c_1)) = g_2(\jt(c_2))$.

Denote the set of edges at $c_1$ in $\tiT_1$ by $\ticE_{c_1}$. Let $g_1^{-1}(\jt(\nu)) = \{\jt(c_1)^i\}_{i=0}^{d-1} \subset \ticE_{c_1}$, written in the cyclic order given by the orientation at $c_1$ with $\jt(c_1)^0 = \jt(c_1)$. Partition $\ticE_{c_1}$ into $d$ pairwise unlinked sets $\ticE_{c_1}^0, \ldots, \ticE_{c_1}^{d-1}$, such that for $0 \leq i <d$, $g_1$ maps $\ticE_{c_1}^i$ bijectively onto $\cE_\nu$, and $\jt(c_1)^i$ is the minimal edge for the cyclic order on $\ticE_{c_1}^i$. For $0 \leq j < d$, define $T_1^j$ as the connected component of
$$
\tiT_1 \setminus \bigcup_{i\neq j} \ticE_{c_1}^i
$$
containing $\jt(c_1)^j$. Observe that $g_1^j := g_1|_{T_1^j} : T_1^j \to T$ is a  tree cover such that $g_1^0$ has exactly $n$ critical vertices ($c_1$ is no longer critical), and $g_1^j$ for $j \neq 0$ is a tree isomorphism. Let $g_2^j := g_2|_{T_2^j} : T_2^j \to T$ be corresponding  tree covers coming from $g_2$. By the induction hypothesis, $\phi$ extends uniquely as a tree isomorphism between $T_1^i$ and $T_2^i$ for $0 \leq i < d$.
\end{proof}

\begin{defn}\label{defn.comb map}
A {\it tree branched cover $h \colon T \to T$} is a tree homomorphism that satisfies the following properties.
\begin{enumerate}
\item For $v \in \cV$, there exists $\deg(v) \geq 1$, referred to as the {\it degree of $v$}, such that every $e \in \cE(h(v))$ has exactly $\deg(v)$ preimages in $\cE(v)$.
\item The set $\cV^{\crit} = \cV^{\crit}(h) \subset \cV$ of vertices whose degrees are at least $2$ is finite.
\item There exists a finite subset $\cE^{\crit}=\cE^{\crit}(h) \subset \cE$ such that for  $e = (v, \tiv) \in \cE$, we have $h(v) = h(\tiv)$ if $e \in \cE^{\crit}$, and $h(v) \neq h(\tiv)$ if otherwise.
\item If $e = (v,v) \in \cE_v$ for some $v \in \cV$, then there exists $e_c \in \cE^{\crit}$ and $n \geq 1$ such that $e = h^n(e_c)$.
\item For $e \in \cE$, there exists $n \geq 0$ such that $h^n(e) \in \cE^{\crit}$.
\end{enumerate}
A vertex $v_c \in \cV^{\crit}$ and an edge $e_c \in \cE^{\crit}$ are said to be {\it critical}. Define $T^{\crit} = T^{\crit}(h)$ as the minimal forward invariant subtree of $T$ containing $\cV^{\crit}$ and $\cE^{\crit}$ (by convention, we set $T^{\crit} := \{v^{\rt}\}$ if $\cV^{\crit}, \cE^{\crit} = \varnothing$).
\end{defn}

\begin{lem}\label{map to root}
Let $h \colon T \to T$ be a tree branched cover. Then $\het(h^n(v))$ is non-increasing with $n \geq 0$, and is eventually $0$ (i.e. $v$ eventually maps to $v^{\rt}$). Consequently, if $\cE^{\crit}_{v^{\rt}} =\varnothing$, then $\cV = \{v^{\rt}\}$ and $\cE = \varnothing$.
\end{lem}

\begin{proof}
Let $v \in \cV$. Observe that $|\Gamma_T(h(v))| \leq |h(\Gamma_T(v))| \leq |\Gamma_T(v)|$. It follows that $\het(h^n(v))$ is non-increasing.

If $\het(v) = 1$, then $\jt(v) \in \cE_{v^{\rt}}$. Thus, $h^n(\jt(v)) \in \cE_{v^{\rt}}$ for $n \geq 0$, and there exists $m \geq 0$ such that $h^m(\jt(v)) \in \cE^{\crit}_{v^{\rt}}$. Then $h^{m+1}(v) = v^{\rt}$.

Proceeding by induction, suppose that for all $u \in \cV$ with $\het(u) < \het(v)$, there exists $n \in \bbN$ such that $\het(h^n(u)) < \het(u)$. By the induction hypothesis, we see that $|h^n(\Gamma_T(v))| < |\Gamma_T(v)|$ for some $n \in N$. Thus, $v$ eventually maps to $v^{\rt}$.
\end{proof}

\begin{defn}
Let $h \colon T \to T$ be a tree branched cover. The {\it generation $\gen(v) \geq 0$ of $v \in \cV$} is the smallest number $n = \gen(v)$ such that $h^n(v) = v^{\rt}$.
\end{defn}

\begin{prop}\label{map ext}
Let $h : T \to T$ and $\tih : \tiT \to \tiT$ be tree branched covers defined on $T = (\cV, \cE, V^{\rt})$ and $\tiT = (\ticV, \ticE, \tiV^{\rt})$ respectively. If there exists a bijection
$$
\phi_{\crit} : T^{\crit}(h) \cup \cE_{T^{\crit}(h)} \to \tiT^{\crit}(\tih) \cup \ticE_{\tiT^{\crit}(\tih)}
$$
that restricts to a tree isomorphism from $T^{\crit}(h)$ to $\tiT^{\crit}(\tih)$ such that $\phi_{\crit} \circ h = \tih \circ \phi_{\crit}$, then $\phi_{\crit}$ extends uniquely to a tree isomorphism $\phi : T \to \tiT$ such that $\phi \circ h = \tih \circ \phi.$
\end{prop}

\begin{proof}
If $\cE^{\crit}_{v^{\rt}} = \varnothing$, the result is trivial by \lemref{map to root}. Proceeding by induction, suppose that the result is true for tree branched covers whose total number of critical vertices and edges are less than $n \geq 1$. Let $c$ be a critical vertex or edge of $h$ with the largest height.

Suppose $c \in \cE^{\crit}$. Denote
$$
\cO := \bigcup_{i=0}^{\infty} h^{-i}(c).
$$
For $e \in \cO$, write $e = (u_e, v_e)$ so that $e$ is the joint of $v_e$, and let $T_e$ be the connected component of $T \setminus \{e\}$ containing $v_e$. Define $T'$ as the connected component of
$$
T \setminus \bigcup_{e \in \cO} T_e
$$
containing $v^{\rt}$. Then $h$ restricted to $T'$ is a Fatou tree branched cover with one less critical edge.

Let $\tic := \phi_{\crit}(c)$. Proceed with the same construction as above, labeling the corresponding objects for $\tih$ by tilde. We again conclude that $\tih$ restricted to $\tiT'$ is a Fatou tree branched cover with one less critical edge. By the inductive hypothesis, $\phi_{\crit}$ extends uniquely to a tree isomorphism $\phi : T' \to \tiT'$ that conjugates $h$ and $\tih$. 

For $e = (u_e, v_e) \in \cO$, define $T_e'$ as the connected component of
$$
T_e \setminus \bigcup_{d \in \cO \cap T_e} T_d
$$
containing $v_e$. If $h^n(e) = c$, then $h^{n+1}|_{T_e' \cup \{e\}}$ is a tree isomorphism from $T_e' \cup \{e\}$ to $T'$. If $\phi$ is already defined on $u_e$, then there exists a unique tree $\tiT_{\tie}'$ attached via the edge $\tie := \phi(e) \in \ticO$ to $\tiu_{\tie} := \phi(u_e)$ such that $\tih^{n+1}|_{\tiT_{\tie}' \cup \{\tie\}}$ is a tree isomorphism from $\tiT_{\tie}' \cup \{\tie\}$ to $\tiT'$. In this case, $\phi$ extends uniquely to a tree isomorphism from $T_e'$ to $\tiT_{\tie}'$ by \propref{cov ext} (uniqueness is guaranteed by the fact that $\phi$ must send $e$ to $\tie$). Since
$$
T = T' \cup \bigcup_{e \in \cO} T_e',
$$
we see that $\phi$ extends uniquely to a tree isomorphism from $T$ to $\tiT$ so that it conjugates $h$ and $\tih$.

The argument in the case $c \in \cV^{\crit}$ is analogous, except the modification to disconnect copies of Fatou trees attached at iterated preimages of $c$ is done as in the proof of \propref{cov ext}.
\end{proof}


\section{Fatou Trees}
\label{sec:fatoutrees}

\subsection{Topological structure on oriented  trees}

By Carath\'eodory, any full compact subset of $\bbC$ with locally connected boundary can be continuously parameterized by $S^1$, and hence is orientable. Henceforth, we endow all such objects with counter-clockwise orientation.

\begin{lem}\label{cover full comp}
Let $V \subset \bbC$ be a full compact set with $\Int(V) \neq \varnothing$, and let $\eta : V' \to V$ be a proper finite-to-one branched cover from a set $V' \subset \bbC$ with all branching points contained in $\Int(V')$. Then $V'$ is also full and compact.
\end{lem}

\begin{proof}
Let $\{c_i'\}_{i=1}^d$ be the set of branching points of $\eta$. For $1 \leq i \leq d$, let $c_i := \eta(c_i')$, and $B_i', B_i$ be the connected components of $\Int(V')$ and $\Int(V)$ containing $c_i'$ and $c_i$ respectively. Then $\overline{B_i}$ and $\overline{B_i'}$ must be closed disks whose boundaries are Jordan curves.

Let $\{C_i\}_{i=1}^k$ and $\{C_i'\}_{i=1}^{k'}$ be the set of closures of the connected components of $V \setminus \bigcup_{i=1}^d B_i$ and $V' \setminus \bigcup_{i=1}^d B_i'$ respectively. Note that each $C_i'$ is mapped homeomorphically to $C_j$ for some $1 \leq j \leq k$. Hence, $C_i'$ is full and compact.

We claim that for $1 \leq i \leq d$ and $1 \leq j \leq k$, $C_i \cap \overline{B_i}$ is either the empty set or a singleton. Otherwise, there is a Jordan arc $\gamma$ contained in $\partial C_i$ whose endpoints are two distinct points in $\partial B_i$. Since $V$ is full, the bounded connected component $K$ of $\bbC \setminus (\gamma \cup \overline{B_i})$ must be contained in $V$, which contradicts the fact that $B_i$ is a connected component of $\Int(V)$.

Since $\eta$ is locally homeomorphic, it follows that for $1 \leq i \leq d$ and $1 \leq j \leq k$, $C_i' \cap \overline{B_i'}$ is either the empty set or a singleton. The fact that $V'$ is full and compact is now immediate.
\end{proof}

\begin{defn}
A {\it Fatou tree} $\cT = (\cV, \cE, V^{\rt})$ is a subset of $\bbC$ with an oriented tree structure that satisfies the following properties.
\begin{enumerate}
\item Every $V \in \cV$ is a full compact subset of $\bbC$ with locally-connected boundary and non-empty interior.
\item If $V, \tiV \in \cV$, $V \neq \tiV$ and $(V, \tiV) \notin \cE$, then $V \cap \tiV = \varnothing$.
\item If $e = (V, \tiV) \in E$ and $V \neq \tiV$, then $e$ is the point in $\bbC$ given by
$$
\{e\} = V \cap \tiV.
$$
\end{enumerate}
\end{defn}

Consider an infinite rooted tree path $\Gamma \sim \{V^i\}_{i=0}^\infty$ in a Fatou tree $\cT$. Its accumulation set in $\bbC$ is given by
$$
\omega(\Gamma) := \bigcap_{n=0}^\infty \cl\left(\bigcup_{i=n}^\infty V^i\right).
$$
If $\omega(\Gamma)= \{x\}$, then we say that $\Gamma$ {\it lands at $x$}. Otherwise, $\omega(\Gamma)$ is referred to as a {\it landing residue}.

For $n \in \bbN$, consider an infinite rooted tree path $\Gamma_n \sim \{V^i_n\}_{i=0}^\infty$. Suppose there exists $K \geq 0$ such that that $V^i_n = V^i_m =: V^i$ for all $n, m$ and $i \leq K$, and $\{V^{K+1}_n\}_{n=1}^\infty$ are pairwise distinct. If
$$
\chi = \bigcap_{N=1}^\infty \cl\left(\bigcup_{n=N}^\infty\bigcup_{i=K+1}^\infty V^i_n\right)
$$
is a non-singleton set, then it is referred to as a {\it ghost limb (at $V^K$)}.

We say that a Fatou tree is {\it starlike} if there are no landing residues or ghost limbs, and every pair of distinct infinite tree paths land at distinct points.

\begin{lem}\label{land extend}
Let $\cT'$ be a subtree of a starlike Fatou tree $\cT$. Then $\partial \cT'$ is locally connected, and $\overline{\cT'}$ is a full compact subset of $\bbC$.
\end{lem}

\begin{defn}
A {\it  Fatou tree cover} $\gamma : \cT \to \ticT$ is a continuous tree cover defined on a Fatou tree $\cT = (\cV, \cE, V^{\rt})$ such that for $V \in \cV$, $\gamma|_V$ is an orientation-preserving $\deg(V)$-to-$1$ proper branched cover from $V$ to $\gamma(V)$, with all branching points (if any) contained in $\Int(V)$. If $\gamma$ is a homeomorphism, it is called a {\it Fatou tree isomorphism}.
\end{defn}

\begin{lem}\label{isom extend}
Let $\cT, \cT'$ be starlike Fatou trees. Then any Fatou tree cover $\gamma : \cT \to \cT'$ has a unique continuous extension $\gamma: \overline{\cT} \to \overline{\cT'}$. Moreover, if $\gamma|_{\cT}$ is a homeomorphism, then so is $\gamma|_{\overline{\cT}}$.
\end{lem}

\begin{defn}\label{defn.fat map}
A {\it Fatou tree branched cover $\eta : \cT \to \cT$} is a continuous tree branched cover that satisfies the following properties.
\begin{enumerate}
\item For $V \in \cV$, $\eta|_V$ is an orientation-preserving $\deg(V)$-to-$1$ proper branched cover from $V$ to $\eta(V)$, with all branching points (if any) contained in $\Int(V)$
\item If $e \in \cE \setminus \cE^{\crit}$, then for any $\tie \in \cE$ with $e \neq \tie$, the edges $e$ and $\tie$ are distinct as points in $\bbC$.
\end{enumerate}
In (2), we allow for the possibility that multiple edges $\{e_i = (V, V_i)\}_{i=1}^{d-1} \subset \cE^{\crit}(V)$ are identified as a single critical point $c \in \partial V$ of degree $\deg(c):= d\geq 2$.
\end{defn}

\subsection{Grand Fatou tree}

Let $\Psi : \bbD \to \cT_0$ be an orientation-preserving homeomorphism, and let $\eta_0 : \cT_0 \to \cT_0$ be given by either
$$
\Psi^{-1}\circ \eta_0\circ \Psi(z) := e^{2\pi i \rho}z
$$
for some $\rho \in (\bbR\setminus \bbQ)/\bbZ$, or
$$
\Psi^{-1}\circ \eta_0\circ \Psi(z) := z^{D_0},
$$
for some $D_0 \geq 2$. The homeomorphism $\Psi$ is referred to as the {\it root uniformization of $\eta$}. For $1 \leq n \leq N$ with $N \in \bbN$, let $\eta_n : \cT_n \to \cT_n$ be a Fatou tree branched cover defined on $\cT_n = (\cV_n, \cE_n, V^{\rt}_n)$ such that $\cT_n$ is starlike, and $V^{\rt}_{n+1} = \overline{\cT_n}$. Then we call $\cT := \cT_N$ a {\it grand Fatou tree}, and $\eta := \eta_N$ {\it a grand Fatou tree branched cover}. The sequence $\{\eta_n\}_{n=0}^N$ is referred to as the {\it Fatou tree tower of $\eta : \cT \to \cT$}.

Denote $\cW_N := \cV_N$. The root $V^{\rt}_N \in \cW_N$ of scale $N$ has a tree structure of scale $N-1$:
$$
(\cV_{N-1}, \cE_{N-1}, V^{\rt}_{N-1})
$$
through the identification $V^{\rt}_N = \overline{\cT_{N-1}}$. This structure can then be transported by the dynamics to any vertex $W_N \in \mathcal{V_N}$ of scale $N$ and produce its tree structure of scale $N-1$:
$$
(\cV_{N-1}(W_N), \cE_{N-1}(W_N), V^{\rt}_{N-1}(W_N)).
$$
The set of vertices of the grand Fatou tree $\cT$ of scale $N-1$ is defined as: 
$$
\cW_{N-1} := \bigcup_{W_N \in \cV_N} \cV_{N-1}(W_N).
$$
Proceeding inductively, we define the vertex set $\cW_n$ of scale $1 \leq n < N$ of $\cT$, and each vertex $W_n \in \cW_n$ has a tree structure of scale $n-1$. See Figures \ref{Fig:1}, \ref{Fig:2} and \ref{Fig:3}. The details of this construction is given below.

\begin{figure}[ht]
\includegraphics[scale=0.5]{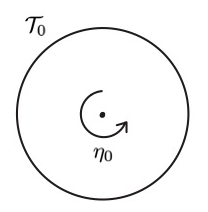}
\caption{A Fatou tree $\cT_0 = (\cV_0, \cE_0, V^{\rt}_0)$ of scale $0$. Note that $\cV_0 = \{V^{\rt}_0\}$ and $\cE_0 = \varnothing$.}
\label{Fig:1}
\end{figure}

\begin{figure}[ht]
\includegraphics[scale=0.5]{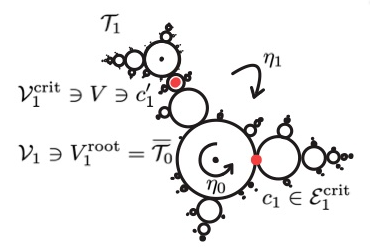}
\caption{A Fatou tree $\cT_1 = (\cV_1, \cE_1, V^{\rt}_1)$ of scale $1$. Note that $V^{\rt}_1 = \overline{\cT_0}$. There are one critical edge $c_1 \in \cE_1^{\crit} \subset \cE_1$ and one critical vertex $V \in \cV_1^{\crit} \subset \cV_1$.}
\label{Fig:2}
\end{figure}

\begin{figure}[ht]
\includegraphics[scale=0.46]{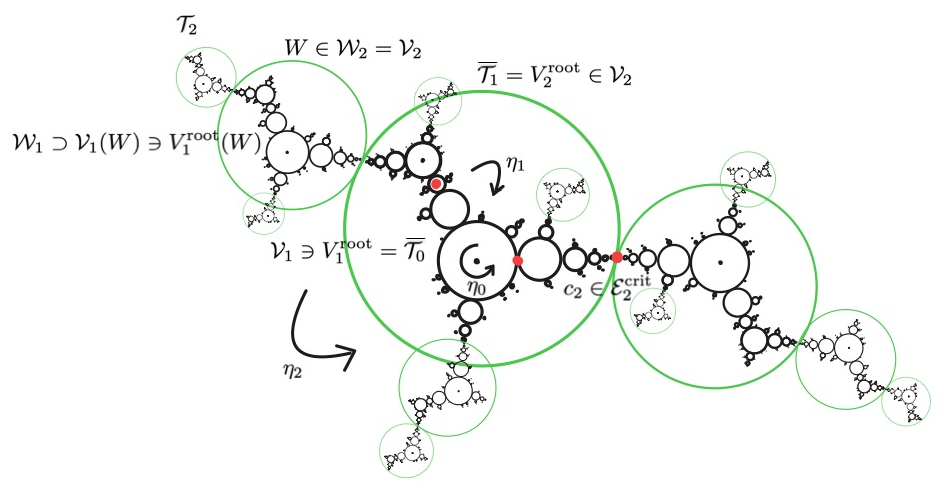}
\caption{A Fatou tree $\cT_2 = (\cV_2, \cE_2, V^{\rt}_2)$ of scale $2$. Note that $V^{\rt}_2 = \overline{\cT_1}$. There is one critical edge $c_2 \in \cE_2^{\crit} \subset \cE_2$. Each vertex $W \in \cV_2$ has a tree structure $(\cV_1(W), \cE_1(W), V^{\rt}_1(W))$ of scale $1$. The collection of all vertices of scale $1$ is $\cW_1$.}
\label{Fig:3}
\end{figure}

Suppose that $\cW_{n+1}$ is a well-defined set such that for $W \in \cW_{n+1}$, there exists a minimal number $\gen(W) \geq 0$ such that $\eta^{\gen(W)}$ is a proper branched cover of $W$ onto $V^{\rt}_{n+1} = \overline{\cT_n}$ with all branch points contained in $\Int(W)$. Through $\eta^{\gen(W)}$, $W$ inherits the graph structure of $\cT_n$. Let $\cV_n(W)$ and $\cE_n(W)$ be the pullback of the vertex set $\cV_n$ and the edge set $\cE_n$ of $\cT_n$ via $\eta^{\gen(W)}$. By \lemref{cover full comp}, each vertex $U \in \cV_n(W)$ is full and compact. Moreover, we have $U_1 := \eta^{\gen(W)}(U) \in \cV_n$. There exists $\gen(U_1) \geq 0$ such that
$$
\eta_n^{\gen(U_1)}(U_1) = V^{\rt}_n = \overline{\cT_{n-1}}.
$$
Hence, $\gen(U) := \gen(W) + \gen(U_1)$. Define
$$
\cW_n := \bigcup_{W \in \cW_{n+1}} \cV_n(W).
$$

To identify $W \in \cW_{n+1}$ as a tree, we must assign one of the vertices in $\cV_n(W)$ as the {\it root} of $W$ (denoted by $V^{\rt}_n(W)$). We proceed inductively. Suppose that $V^{\rt}_{n+1}(W')$ and $V^{\rt}_n(W_1)$ are defined for all $W' \in \cW_{n+2}$ and $W_1 \in \cW_{n+1}$ with $\gen(W_1) < \gen(W)$. Let $W_1 := \eta(W) \in \cW_{n+1}$, and let $W', W_1' \in \cW_{n+2}$ be the unique trees such that $W \in \cV_{n+1}(W')$ and $W_1 \in \cV_{n+1}(W_1')$. By the induction hypothesis, the roots $V^{\rt}_{n+1}(W')$ and $V^{\rt}_{n+1}(W_1') = \eta(V^{\rt}_{n+1}(W'))$ are well-defined. Hence, so are the joints $\jt_1(W)$ and $\jt_1(W_1) = \eta(\jt_1(W)) \in \partial W_1$. If $\deg(\eta|_W) =1$, then define
$$
V^{\rt}_n(W) := (\eta|_W)^{-1}(V^{\rt}_n(W_1)).
$$

Suppose instead that $\deg(\eta|_W) \geq 2$, so that $(\eta|_W)^{-1}(V^{\rt}_n(W_1))$ has (possibly) multiple connected components $\{D_i\}_{i=1}^d$. Since $W$ is full and compact, there is an orientation-preserving Carath\'eodory parameterization $\gamma : \bbR/\bbZ \to \partial W$ of its boundary. Denote
$$
X := \gamma^{-1}(\jt_1(W))
\matsp{and}
Y_i := (\gamma)^{-1}(\partial D_i).
$$
Note that $Y_i$'s are compact and pairwise disjoint sets. Moreover, either $X \cap Y_i = \varnothing$ or $X \subset Y_i$. If $X$ is a singleton set $\{x\}$ , then define $V^{\rt}_n(W) := D_j$, where $1 \leq j\leq d$ is the unique number such that $\gamma(x+t) \in Y_j$ with minimal $t \geq 0$.

Suppose instead that $|X| > 1$. Then there exists a rooted tree path $\Gamma_1 \sim \{U_i\}_{i=0}^{K+1} \subset W_1$ with $U_0 = V^{\rt}_n(W_1)$ and $U_K \cap U_{K+1} = \{\jt_1(W_1)\}$. Let $C_1^{\rt}$ be the union of $U_K$ with the connected component of $W_1 \setminus U_K$ containing $V^{\rt}_n(W_1)$, and let $C^{\rt}$ be the connected component of $(\eta|_W)^{-1}(C_1^{\rt})$ containing $\jt_1(W)$. Observe that $C^{\rt}$ is a full compact subset of $\bbC$, and that $\partial C^{\rt} \setminus \{\jt_1(W)\}$ is connected. Thus, for the Carath\'eodory parametrization $\gamma_{C^{\rt}} : \bbR/\bbZ \to \partial C^{\rt}$, $\gamma_{C^{\rt}}^{-1}(\jt_1(W))$ is a singleton set $\{x\}$. We define $V^{\rt}_n(W) := D_j \subset C^{\rt}$, where $1 \leq j\leq d$ is the unique number such that $\gamma_{C^{\rt}}(x+t) \in \gamma_{C^{\rt}}(D_j)$ with minimal $t \geq 0$.

Let $0 \leq n < N$. A {\it nested tree sequence $\cS = \{W_i\}_{i=n}^N$} is a sequence such that $W_N = \cT$, and for $n \leq i < N$, we have $W_i \in \cW_i$ and $W_i \subset W_{i+1}$. Note that for $W \in \cW_n$, there exists a unique nested tree sequence $\cS(W) = \{W_i\}_{i=n}^N$ such that $W_n = W$. We denote $\cS_m(W) := W_m$ for $n \leq m \leq N$.

Let $W \in \cW_0 \setminus V^{\rt}_0$. Denote $W_i := \cS_i(W)$ for $0 \leq i \leq N$. If $W \neq V^{\rt}_0(W_1)$, then the joint $\jt(W)$ of $W$ in $W_1$ is well-defined. Otherwise, let $1 < k \leq N$ be the smallest number such that $W_k \neq V^{\rt}_k(W_{k+1})$. Define $\jt_{k+1}(W_k) \in \partial W_k$ as the joint of $W_k$ in $W_{k+1}$. Proceeding inductively, if $\jt_{n+1}(W_n) \in \partial W_n$ is defined, let $\jt_n(W_{n-1}) \in \partial W_{n-1}$ be the unique edge contained in the rooted tree path $\Gamma_{W_n}(\jt_{n+1}(W_n))$. Define $\jt(W) := \jt_0(W_0)$.

\subsection{Combinatorial rigidity of Fatou trees}

The {\it root parameterization $\psi_{\cT_0} : \bbR/\bbZ \to \partial \cT_0$ of $\eta$} is defined by $\psi_{\cT_0} := \Psi|_{\partial \bbD}$ (where $\partial \bbD$ is identified with $\bbR/\bbZ$). To parameterize $\partial V$ for $V \in \cW_0$, we use the following elementary result. 

\begin{lem}
Let $V \in \cW_0$. Denote the degree of $\eta^{\gen(V)}|_V$ by $d_V \geq 1$. Then there exists a unique homeomorphism $\psi_V : \partial \bbD \to \partial V$ such that
$$
\psi_{\cT_0}^{-1} \circ \eta^{\gen(V)}\circ \psi_V(t) = d_Vt;
$$
and $\jt(V) \in \psi_V([0, 1/d_V))$, where $\jt(V)$ is the joint of $V$.
\end{lem}

For $x \in \cT$, we define its {\it combinatorial address $\adr(x)$} as follows. Consider a nested tree sequence $\cS = \{W_i\}_{i=0}^N$. First, define
$$
\adr_0(x) := \psi_{W_0}^{-1}(x)
\matsp{for}
x \in \partial W_0.
$$
Proceeding inductively, for $0 \leq n < N$, suppose that $\adr_n(x)$ defined for all $x \in \partial W_n$. If $x \in \partial W_{n+1}$, then either $x \in \partial V$ for some $V \in \cV_n(W_{n+1})$ or $x$ is the landing point of some infinite rooted tree path in $W_{n+1}$.

Let $\Gamma \sim \{V^i\}_{i=0}^K$ be a rooted tree path in $T_{n+1}$ of length $0 \leq K \leq \infty$ with the joint set $\{e^i\}_{i=0}^{K-1}$. Define
$$
[\Gamma]_{n+1} = [V^K]_{n+1} := (\adr_n(e^0), \adr_n(e^1),\ldots, \adr_n(e^{K-1}))
$$
if $K < \infty$, and
$$
[\Gamma]_{n+1} = \adr_{n+1}(\omega(\Gamma)) := (\adr_n(e^0), \adr_n(e^1), \ldots).
$$
if $K = \infty$. If $K < \infty$, then for $x \in \partial V^K$, define
$$
\adr_{n+1}(x) := ([V^K]_{n+1}; \adr_n(x)).
$$

Let $x \in \cT$. If $x \in \partial \cT$, define $\adr(x) := \adr_N(x)$. If $x \in \Int V$ for some $V \in \cW_0$, let $\cS(V) = \{W_i\}_{i=0}^N$ be the nested tree sequence at $V = W_0$, and define
$$
\adr(x) := ([W_{N-1}]_N ; \ldots; [W_1]_2 ; [W_0]_1).
$$

\begin{defn}[Critical address]
\label{Def:CritAdr}
The {\it critical address of $\eta : \cT \to \cT$ (with respect to $\psi$)} is
$$
\adr_{\crit}(\eta) := \{\adr(c) \; | \; c \in \cT \text{ is a critical point of } \eta\}.
$$
\end{defn}

\begin{thm}\label{tree rigid}
Let $\eta : \cT \to \cT$ and $\tieta : \ticT \to \ticT$ be grand Fatou tree branched covers with root uniformizations $\psi_{\cT_0} : \bbR/\bbZ \to \partial \cT_0$ and $\tipsi_{\ticT_0}: \bbR/\bbZ \to \partial \ticT_0$ respectively. Let $\Phi_0 := \tipsi_{\ticT_0} \circ \psi_{\cT_0}^{-1}$. If $\Phi_0$ conjugates $\eta|_{\partial \cT_0}$ with $\tieta|_{\partial \ticT_0}$, and $\adr_{\crit}(\eta) = \adr_{\crit}(\tieta)$, then $\Phi_0$ extends uniquely to a homeomorphism $\Phi : \partial \cT \to \partial \ticT$ that conjugates $\eta|_{\partial \cT}$ with $\tieta|_{\partial \ticT}$.
\end{thm}

\begin{proof}
Throughout the proof, objects corresponding to $\tieta$ will be denoted with a tilde.

Let $\{\eta_n : \cT_n \to \cT_n\}_{n=0}^N$ be the Fatou tree tower of $\eta$. Consider a sequence of Fatou trees $\{T_i\}_{i=0}^k$ with pairwise disjoint interiors and a decreasing sequence
$$
N \geq n_0 > n_1 > \ldots > n_k \geq 1
$$
such that $T_0$ is a subtree of scale $n_0$ of
$$
W_0 := V^{\rt}_{n_0+1} = \overline{\cT_{n_0}}
$$
and the following properties hold for $1 \leq i \leq k$:
\begin{itemize}
\item There exists $W_i \in \cW_{n_i+1}$ such that $T_i$ is a subtree of scale $n_i$ of $W_i$.
\item There exists a vertex $W_{i+1}' \in \cV_{n_i}(W_i)$ of $W_i$ such that $W_{i+1}' \supseteq W_{i+1}$ and $W_{i+1}'$ is attached to $T_i$ by an edge.
\end{itemize}
In particular, $W_1' \in \cV_{n_0}= \cV_{n_0}(\cT_{n_0})$. Proceeding by induction, suppose that
\begin{itemize}
\item $\Phi$ extends uniquely as a conjugacy to the boundary of the subtree $\cT_{n_0}' \supset T_0$ of $\cT_{n_0}$ containing all vertices of generation less than $\gen(W_1')$; and
\item $\Phi$ extends uniquely to $\bigcup_{i=0}^k \partial T_i$.
\end{itemize}
If $\cT_{n_0}'$ contains $\cT_{n_0}^{\crit}$, then $\ticT_{n_0}'$ contains $\ticT_{n_0}^{\crit}$ since $\adr_{\crit}(\eta) = \adr_{\crit}(\tieta)$. Thus, by Propositions \ref{map ext} and \ref{isom extend}, $\Phi$ extends uniquely to $\partial \cT_{n_0}$ so that
$$
\Phi \circ \eta|_{\partial \cT_{n_0}} = \tieta|_{\partial \ticT_{n_0}}\circ \Phi.
$$
In this case, reapply the inductive procedure after replacing $\{T_i\}_{i=0}^k$ with the singleton set $\{\hat{T}_0 = V^{\rt}_{n_0+1}\}.$

Suppose this is not the case. For $1 \leq i \leq k$, let $W_i^{\crit}$ be the minimal subtree of $W_i$ containing all the critical vertices of the Fatou tree cover $\eta|_{W_i}$. If $T_i$ contains $W_i^{\crit}$, then $\tiT_i$ contains $\tiW_i^{\crit}$ since $\adr_{\crit}(\eta) = \adr_{\crit}(\tieta)$. Thus, by Propositions \ref{cov ext} and \ref{isom extend}, $\Phi$ extends uniquely to $\partial W_i$ so that
$$
\Phi \circ \eta|_{\partial W_i} = \tieta|_{\partial \tiW_i} \circ \Phi.
$$
In this case, let $0 \leq j \leq k$ be the smallest number such that $T_j \supset W_j^{\crit}$. If $n_j +1 < n_{j-1}$, then reapply the inductive procedure after replacing $\{T_i\}_{i=0}^k$ by $\{\hat{T}_i\}_{i=0}^j$, where $\hat{T}_i = T_i$ for $i < j$ and $\hat{T}_j = W_j \in \cW_{n_j+1}$. If $n_j +1 = n_{j-1}$, then note that $W_j = W_j' \in \cV_{n_j}(W_{j-1})$. In this case, reapply the inductive procedure after replacing $\{T_i\}_{i=0}^k$ by $\{\hat{T}_i\}_{i=0}^{j-1}$, where $\hat{T}_i = T_i$ for $i < j-1$ and $\hat{T}_{j-1} = T_{j-1} \cup W_j$.

Suppose this is not the case. Consider $W_{k+1}' \in \cV_{n_k}(W_k)$ attached to $T_k$ by an edge. Let $T_{k+1} := V^{\rt}_1(W_{k+1}')$, and extend $\Phi$ uniquely to $T_{k+1}$ by
$$
\Phi|_{T_{k+1}} := \psi_{\tiT_{k+1}}\circ \psi_{T_{k+1}}^{-1}.
$$
Since $\adr_{\crit}(\eta) = \adr_{\crit}(\tieta)$, it follows that
$$
\Phi \circ \eta|_{\partial T_{k+1}} = \tieta|_{\partial \tiT_{k+1}} \circ \Phi.
$$
If $n_k >1$, then reapply the inductive procedure after replacing $\{T_i\}_{i=0}^k$ by $\{T_i\}_{i=0}^{k+1}$. If $n_k = 1$, then note that $T_{k+1} = W_{k+1}' \in \cV_1(W_k)$. In this case, reapply the inductive procedure after replacing $\{T_i\}_{i=0}^k$ by $\{\hat T_i\}_{i=0}^k$, where $\hat T_i = T_i$ for $i < k$ and $\hat T_k = T_k \cup T_{k+1}$.

Since the sets of critical vertices and edges of $\eta$ and $\tieta$ are finite, this process eventually extends $\Phi$ uniquely, in finitely many iterates, to a homeomorphism from $\partial \cT$ to $\partial \ticT$ that conjugates $\eta$ and $\tieta$.
\end{proof}


\section{Topological Models of Atomic Siegel Julia Sets}
\label{sec:topmodel}

Let $f : \bbC \to \bbC$ be an atomic Siegel polynomial of bounded type with rotation number $\rho \in (\bbR\setminus \bbQ)/\bbZ$. Recall that this means $f$ satisfies the following properties:
\begin{enumerate}
\item $f$ has a fixed Siegel disk $\Delta_f \subset \bbC$ of rotation number $\rho$.
\item $f$ has a trivial lamination of rational external rays.
\end{enumerate}
In this section, we show that the filled Julia set $K_f$ of such a polynomial is a grand Fatou tree, and that $f$ is a grand Fatou tree cover with the Fatou tree tower $\{f_n : \cT_n \to \cT_n\}_{n=0}^N$.

Let $\cT_0 := \Delta_f$. The root parameterization of $f$ is given by the quasisymmetric map $\psi_{\Delta_f} : \bbR/\bbZ \to \partial \Delta_f$ that linearizes $f|_{\partial \Delta_f}$ to rotation by $\rho$. Proceeding by induction, suppose that for $n \geq 1$, $f_{n-1} : \cT_{n-1} \to \cT_{n-1}$ is a well-defined Fatou tree branched covering such that $\cT_{n-1}$ is a starlike, invariant subset of $K_f$.

Set $V^{\rt}_n := \overline{\cT_{n-1}}$. Let $\cW_n$ be the set of compact subsets of $K_f$ such that if $W \in \cW_n$, then there exists the smallest number $m:= \gen(W) \geq 0$ such that $f^m|_W$ is a proper branched covering map from $W$ to $V^{\rt}_n$ with all branching points contained in $\Int(W)$.

\begin{prop}\label{f edge}
Let $W, \tiW \in \cW_n$. Then either $W \cap \tiW = \varnothing$ or $W \cap \tiW = \{e\}$, where $e$ is an iterated preimage of a critical point. Consequently, if $\partial V^{\rt}_n$ contains no critical points of $f$, then $V^{\rt}_n = K_f$.
\end{prop}

\begin{proof}
By \lemref{cover full comp}, $W, \tiW \in \cW_n$ are full and compact. Since $K_f$ is full and compact, the fact that $W \cap \tiW = \varnothing$ or $W \cap \tiW = \{e\}$ can be proved the same way as in \lemref{cover full comp}. 

Suppose that $W \cap \tiW = \{e\}$, and that $e_i := f^i(e)$ is not a critical point of $f$ for $i \geq 0$. Then there exists a neighborhood $U_i \ni e_0$ such that $f^i|_{U_i}$ is conformal. Since $f^k(W) = f^k(\tiW) = V_n^{\rt}$ for some $k \geq 0$, this is a contradiction.
\end{proof}

Define $\cT_n$ as the connected component of $\bigcup_{W \in \cW_n} W$ containing $V^{\rt}_n$. Let $\cV_n := \{W \in \cW_n \; | \; W \subset \cT_n\}$ and $\cE_n := \{W \cap \tiW \; | \; W, \tiW \in \cV_n\}$. Lastly, let $f_n := f|_{\cT_n}$.

The following statement is an immediate consequence of \propref{f edge}.

\begin{prop}\label{f is tree}
The set $\cT_n = (\cV_n, \cE_n, V^{\rt}_n)$ is a Fatou tree, and the map $f_n : \cT_n \to \cT_n$ is a Fatou tree branched covering.\qed
\end{prop}

\begin{prop}\label{starlike}
The Fatou tree $\cT_n$ is starlike.
\end{prop}

\begin{proof}
Suppose that $\chi$ is a landing residue of an infinite rooted tree path in $\cT_n$ or a ghost limb at some vertex $V \in \cV_n$. Recall the definition of the initial puzzle $\cZ_0$ given in \eqref{eq:0 puzzle}. Since $f^i(\chi) \cap \cZ_0 = \varnothing$ for all $i \geq 0$, it follows that $\chi \cap \cZ_n =\varnothing$. Hence, $\chi$ is contained in a single fiber. By \thmref{Thm:DynRigidity}, $\chi$ is a singleton.
\end{proof}

By Propositions \ref{f is tree} and \ref{starlike}, we conclude:

\begin{thm}\label{poly is fatou tree}
An atomic Siegel polynomial $f$ restricted to its filled Julia set $K_f$ is a grand Fatou tree branched cover.\qed
\end{thm}

\begin{defn}[Combinatorial equivalence]
\label{Def:CombEq}
Two atomic Siegel polynomials $f, \tif : \bbC \to \bbC$ of bounded type are said to be {\it combinatorially equivalent} if
\begin{enumerate}
    \item 
    $\deg(f) = \deg(\tif\,)$, $\rho(f)=\rho(\tif\,)$; 
    \item 
    $\adr_{\crit}(f)=\adr_{\crit}(\tif)$ (see Definition~\ref{Def:CritAdr});
    \item 
    the orbits of each pair $c \in \Crit(f) \cap \inter(K_f)$, $\tilde c\in \Crit(\tif\,) \cap \inter(K_{\tif})$ of the corresponding critical points land in $\Delta_f$, $\Delta_{\tif}$ in the same conformal position with respect to the uniformizing coordinates. 
\end{enumerate}
\end{defn}

The last condition implies that combinatorially equivalent bounded type Siegel polynomials $f$ and $\tif$ are conformally conjugate restricted to the interior of their Fatou sets. By Theorems \ref{tree rigid} and \ref{poly is fatou tree}, this conjugacy can be extended to a global topological conjugacy.

\begin{thm}[Topological rigidity]\label{poly topo rigid}
Let $f, \tif : \bbC \to \bbC$ be atomic Siegel polynomials of bounded type. If $f$ and $\tif$ are combinatorially equivalent, then they are topologically conjugate, and, restricted to the Fatou sets, this conjugacy is conformal.\qed
\end{thm}

\section{Proof of Main Theorem~\eqref{it:main:3}}
\label{sec:top2qc}

The proof of Main Theorem~\eqref{it:main:3} splits into two steps: 1) to show that combinatorial equivalence implies topological conjugacy (which is conformal in the Fatou set); 2) to show that topological conjugacy implies quasiconformal conjugacy. This, together with Main Theorem~\eqref{it:main:2}, implies Main Theorem~\eqref{it:main:3}. This will conclude the proof.

Step 1) is done in Theorem~\ref{poly topo rigid}. In this section, we carry out step 2) and prove that topological conjugacy implies quasiconformal conjugacy, which must be, in turn, conformal. The main result of this section is the following theorem:

\begin{thm}[Conformal rigidity]
\label{Thm:Rigidity}
Let $f$ and $\tf$ be a pair of atomic Siegel polynomials of bounded type. If they are topologically conjugate and, restricted to the Fatou sets, this conjugacy is conformal, then $f$ and $\tif$ are affinely conjugate.
\end{thm}

Let us denote the topological conjugacy between $f$ and $\tf$ as $\tau$, e.g., $\tf \circ \tau = \tau \circ f$. By assumption, $\tau$ is conformal in the basin of infinity, and hence, coincides with the B\"ottcher coordinate there. Furthermore, again by assumption, $\tau$ is conformal on the grand orbit of the Siegel disk $\Delta_f$.

Furthermore, $\tau \colon \partial \Delta_f \to \partial \Delta_{\tilde f}$ is quasi-symmetric because the maps $f$ and $\tf$ are quasiconformally conjugate to the same Blaschke model. 

We say that a homeomorphism $h \colon P \to \tilde P$ between two puzzle pieces of $f$ and $\tf$ \emph{respects the boundary marking} if it extends continuously to the closures of $P$ and $\tilde P$, and this extension, restricted to $\partial P$, agrees with $\tau$.

We say that two puzzle pieces $P$ of $f$ and $\tilde P$ of $\tilde f$ are \emph{corresponding} if $\tau(P) = \tilde P$. The corresponding puzzle pieces are necessarily of the same depth and have the same combinatorial structure.

\begin{lem}
\label{Lem:StartingQC}
For every pair of corresponding puzzle pieces $P$ and $\tilde P$, there exists a quasiconformal homeomorphism $\phi_P \colon P \to \tilde P$ that respects the boundary marking.
\end{lem}

\begin{proof}
We work in the Blaschke models $\bff, \tilde \bff$ for the polynomials $f$ and $\tilde f$. In this model, the closures of $P, \tilde P$ correspond to a pair of puzzle pieces $\Bp, \tilde{\Bp}$ (see \cite{Y}). In this way, it is enough to establish the claim of the lemma for a pair of the corresponding puzzle pieces $\Bp$ and $\tilde{\Bp}$ with the boundary marking $\boldsymbol{\tau} \colon \partial \Bp \to \partial \tilde{\Bp}$. The lemma will then follow by the quasiconformal equivalence between the Blaschke model and the original coordinates. 

The Siegel part of the boundary $\partial \Bp$ consists of countably many analytic arcs meeting at a definite angle, while the parts of $\partial \Bp$ in the basins of $\infty$ and $0$ consist of finitely many analytic arcs meeting at a definite angle as well.

Let $(b_i) \subset \partial \Bp$ be a collection of Siegel arcs landing at a point $p$; denote by $t_i$ a common point between $b_i$ and $b_{i+1}$. In this way, $\lim_{i \to \infty} t_i = p$. Let $r$ be the arc of an external ray in $\partial \Bp$ landing at $p$.  Assume that $p$ is a repelling fixed point (we will discuss periodic and pre-periodic cases later). Choose a sequence $(s_i) \subset r$ of points linked by the dynamics, i.e., $\bff(s_{i+1}) = s_i$; in this way, $\lim_{i \to \infty} s_i = p$. For a given $i_0$, we connect $s_{i_0}$ to $t_{i_0}$ with an analytic arc $a_{i_0} \subset \inter \Bp$ that makes a definite angle with the arcs $b_{i_0}$ and $r_{i_0}:=[s_{i_0}, s_{i_0+1}] \subset r$. Since the angles between $b_i$ and $b_{i-1}$ are definite, we can choose $a_{i_0}$ so that it also makes a definite angle with $b_{i_0-1}$. Furthermore, choose $i_0$ large enough so that for all $i \ge i_0$, all the arcs $b_i$ and $[s_i, s_{i+1}]$ for $i \ge i_0$, as well as the curve $a_{i_0}$, lie in the domain of the linearizing coordinates of $p$. This allows us to construct a sequence $a_{i_0}, a_{i_0+1}, \ldots$ of arcs by lifting $a_{i_0}$ under the local dynamics near the repelling periodic point. In this way, we obtain a sequence $(R_i)_{i \ge i_0} \subset \Bp$ of closed shrinking curvilinear rectangles bounded by four quasi-symmetric arcs $a_i, b_i, a_{i+1}, r_i$. In the same way, we construct the sequence $(\tilde R_i)_{i \ge i_0} \subset \tilde{\Bp}$ of the corresponding shrinking rectangles bounded by the corresponding quasi-symmetric arcs $\tilde a_i, \tilde b_i, \tilde a_{i+1}, \tilde r_i$. 

By construction, we have a pair quasi-symmetric maps $\boldsymbol{\tau} \colon b_i \to \tilde b_i$ and $\boldsymbol{\tau} \colon r_{i} \to \tilde r_{i}$ that come from the marking $\boldsymbol{\tau}$. Let $h_0 \colon a_{i_0} \to \tilde a_{i_0}$ be any quasi-symmetric homeomorphism, and define inductively $h_i \colon a_{i} \to \tilde a_{i}$ to be the lift of $h_{i-1}$ under the linearized dynamics in a neighborhood of $p$. In this way, for each $i \ge i_0$, we obtained a quasi-symmetric identification $H_i \colon \partial R_i \to \partial \tilde R_i$ that respects the marking on the two opposite sides:
\begin{equation*}
H_i(z) := \left\{
\begin{aligned}
&h_i(z), &\text{if} &\quad z \in a_i\\
&\boldsymbol{\tau}(z), &\text{if} &\quad z \in b_i\\
&h_{i+1}(z), &\text{if} &\quad z \in a_{i+1}\\
&\boldsymbol{\tau}(z), &\text{if} &\quad z \in r_i.\\
\end{aligned}
\right.
\end{equation*}

Extend this quasi-symmetric map to a homeomorphism $H_i \colon R_i \to \tilde R_i$ that is quasiconformal in the interior of the rectangle. For $i \to \infty$, these maps glue nicely on the union $R_p := \bigcup_{i \ge i_0} R_i$ into a homeomorphism 
\[
H_p \colon R_p \to \tilde R_{\tilde p}, \quad H_p|R_i := H_i,
\]
where $\tilde p$ is the point corresponding to $p$. This is a well-defined homeomorphism that is quasiconformal in the interior of $R_p$ and that respects the boundary marking $\boldsymbol{\tau}$ on $\partial R_p \cap \partial \Bp$. 

If $p$ is periodic or pre-periodic, then we construct analogous maps by going to a suitable iterate that fixes $p$ (and possibly shrinking the domain of the linearizing coordinates depending on the depth of $\Bp$). In this way, if $(p_k)$ are all the vertices on $\partial \Bp$ where the Siegel bubble rays meet external rays, we get a collection of homeomorphisms
\[
H_{p_k} \colon R_{p_k} \to \tilde R_{\tilde p_k}
\]
quasiconformal in the interior of each $R_{p_k}$.

The domain $\Bp' := \Bp \sm \bigcup_k \inter(R_{p_k})$ is a finitely-sided closed topological polygon bounded by analytic arcs meeting at definite angles. On each of its sides, we have a quasisymmetric map (that respects the boundary marking on its portion within $\partial \Bp$). We extend this map to a homeomorphism $H_0 \colon \Bp' \to \tilde \Bp'$ quasiconformal in the interior of its domain of definition. 

Finally, we define the map $\Psi \colon \Bp \to \tilde \Bp$ as $H_{p_k}$ in $R_{p_k}$ for every $k$ and as $H_0$ in $\Bp'$. By construction, this is a homeomorphism that is quasiconformal in the interior that respects the boundary marking, as required.
\end{proof}

Using Lemma~\ref{Lem:StartingQC}, we can prove the following proposition. The proof is similar to the proofs of \cite[Proposition 4.1]{KvS} and \cite[Lemma 6.10]{DS}.

\begin{prop}
\label{Prop:CombEquivalent}
If $f$ and $\tf$ are two topologically conjugate, atomic Siegel polynomials of bounded type, then in Theorem~\ref{Thm:Extract} one can choose dynamically natural complex box mappings $F \colon \U \to \V$ and $\tilde F \colon \tilde \U \to \tilde \V$ to be combinatorially equivalent in the sense of Definition~\ref{DefA:CombEquivBoxMappings}, with respect to some quasiconformal homeomorphism $H \colon \V \to \tilde \V$ such that $H(\U) = \tilde \U$ and $\tilde F \circ H = H \circ F$ on $\partial \U$. \qed
\end{prop}

Because the box mappings extracted for Siegel polynomials are non-renormalizable, Theorems~\ref{Thm:BoxLineFields} and~\ref{Thm:BoxParamRigidity} imply that $F$ and $\tilde F$ are quasiconformally conjugate and support no invariant line fields on their non-escaping sets. Let us call this quasiconformal homeomorphism $\Phi$. Note that $\Phi$ agrees with $\tau$ on the boundary of puzzle pieces. 

Now we want to globalize $\Phi$ using the Spreading Principle \cite{CDKvS}. Before stating it, we introduce some notation. For completeness, we will also provide a proof, following \cite{CDKvS}. 

Recall that $\Crit^\star(f)$ is the set of critical points of $f$ in $J_f$ that do not lie in the puzzle graph of any depth and $\Crit^\star(\tf)$ is the corresponding set for $\tf$, and $\fs$ is the smallest depth of the puzzle graph that contains all points in $\left(\Crit(f) \cap J_f\right) \sm \Crit^\star(f)$ (we assume $\fs=0$ if the latter set is empty).

For a finite union $A$ of disjoint puzzle pieces, we define its \emph{depth} as the smallest depth of puzzle pieces in $A$. 

\begin{thm}[Spreading Principle]
\label{Thm:Spreading}
\label{Thm:Spreading_bis}
Let $\Y$ be a nice puzzle neighborhood of $\Crit^\star(f)$ of depth at least $\fs$, and let $\tilde{\Y}$ be the corresponding neighborhood of $\Crit^\star(\tf)$ for $\tf$. Suppose that there exists a $K$-quasiconformal homeomorphism $\phi \colon \Y \to \tilde \Y$ that respects the boundary marking induced by the conjugacy $\tau$. Then $\phi$ extends to a $K$-quasiconformal homeomorphism $\Phi \colon \Cc \to \Cc$ such that:
\begin{enumerate}
\item 
\label{It:SP1}
$\Phi=\phi$ on $\Y$;
\item 
\label{It:SP2}
for each $z\not\in \Y$, $$\tf\circ\Phi(z)=\Phi\circ f(z);$$
\item 
\label{It:SP3}
the dilatation of $\Phi$ vanishes on $\Cc \sm \DomL(\Y)$, where $\DomL(\Y)$ is the union of $\Y$ and all the landing domains to $\Y$;
\item 
\label{It:SP4}
$\Phi(P)=\tilde{P}$ for every puzzle piece $P$ that is not contained in $\DomL(\mathcal Y)$,
					and $\Phi \colon P\rightarrow\tilde{P}$ respects the boundary marking induced by $\tau$.
\end{enumerate}
\end{thm}
\begin{proof}  
Denote by $\Pp_n$ the collection of puzzle pieces of depth $n \ge 0$ for $F$.

Let $\mathcal C$ be a set of all puzzle pieces in $\Pp_0$ as well as all critical puzzle pieces that are not contained in $\Y$ but do intersect $\Crit^\star(f)$. This is a finite set. Using Lemma~\ref{Lem:StartingQC}, for each $C \in \mathcal C$ let us pick a quasiconformal map $h_C \colon C \to \tilde C$ that respects the boundary marking induced by $\tau$; suppose $K_0$ is the maximum over the quasiconformal dilatations of $h_C$ for $C \in \mathcal C$.

Let $Q$ be a component of $\DomL(\Y)$, and $k \ge 0$ be the landing time of the orbit of $Q$ to $\Y$. Since $\Crit^\star(f) \subset \Y$ and the depth of $\Y$ is larger than $\fs$, the map $f^k \colon Q \to f^k(Q)$ is a conformal isomorphism if $k \ge 1$. Therefore, we can pull back $\phi$ and define a $K$-quasiconformal homeomorphism $\phi_Q \colon Q \to \tilde Q$ by the formula $\tf^k \circ \phi_Q = \phi \circ f^k$. Since $\phi$ respects the boundary marking induced by $\tau$, so does $\phi_Q$.

Let us now inductively define a nested sequence of sets $Y_0 \supset Y_1 \supset Y_2 \ldots$ so that:
\begin{itemize}
\item
$Y_0$ is the union of all puzzle pieces in $\Pp_0$;
\item
$Y_{n+1}$ is the subset of $Y_n$ consisting of puzzle pieces of depth $n+1$ that are not components of $\DomL(\Y)$.  
\end{itemize}

For each puzzle piece $Q \in \Pp_n$ that is contained in $Y_n$, it follows that none of the pieces in the orbit $Q, F(Q), \ldots$ can lie in $\Y$ (because otherwise $Q$ would lie in $Y_{\ell-1} \sm Y_\ell$ for some $\ell \le n$). Therefore, there exists the minimal $k = k(Q) \ge 0$ so that $f^k(Q) \in \mathcal C$. By a similar reasoning as above, using the conformal map $f^k \colon Q \to f^k(Q)$ one can pull back and define the homeomorphism $H_Q \colon Q \to \tilde Q$ by the formula $\tf^k \circ H_Q = h_{f^k(Q)} \circ f^k$, where $h_{f^k(Q)}$ was defined above. In this way, $H_Q$ is a $K_0$-quasiconformal map that respects the boundary marking induced by $\tau$.

Now define a sequence of homeomorphisms $(\Phi_n)_{n \ge 0}$, $\Phi_n \colon \Cc \to \Cc$, as follows. Set 
\begin{equation*}
\Phi_{0}(z) = \left\{
\begin{aligned}
&\tau(z) \quad &\text{ if }&z \in \Cc \,\sm Y_0,\\
&h_Q(z) \quad &\text{ if }&z \in Q \in \Pp_0.
\end{aligned}
\right.
\end{equation*}
For each $n \ge 0$, define
\begin{equation*}
\Phi_{n+1}(z) = \left\{
\begin{aligned}
& \Phi_n(z) \quad &\text{ if }&z \in \Cc \,\sm \bigcup_{P \in \Pp_{n}} P,\\
& \tau(z) \quad &\text{ if } &z \in P \,\sm \bigcup_{Q \in \Pp_{n+1}} Q, \text{ where } P \in \Pp_{n}\\
&H_Q(z) \quad &\text{ if }&z \in Q \in \Pp_{n+1} \text{ and }Q \subset Y_{n+1},\\
&\phi_Q(z) \quad &\text{ if }&z \in Q \in \Pp_{n+1}\text{ and }Q \not\subset Y_{n+1}.
\end{aligned}
\right.
\end{equation*}

By the Gluing Lemma (see, e.g., \cite[Lemma 8.2]{CDKvS}), the map $\Phi_n \colon \Cc \to \Cc$ is a $\max\{K, K_0\}$-quasiconformal homeomorphism for each $n \ge 0$. Note that the sequence $(\Phi_n)$ eventually stabilizes on $\Cc \,\sm \bigcap_n Y_n$. The set $\bigcap_n Y_n$ consists of points in $J_f$ which are never mapped into $\Y$, and hence has zero Lebesgue measure (Theorem~\ref{Thm:AccZero}). Hence $\Cc \,\sm \bigcap_n Y_n$ is a dense set. We conclude that the sequence $(\Phi_n)$ converges to the limiting $K$-quasiconformal homeomorphism $\Phi \colon \Cc \to \Cc$. For this homeomorphism the properties \eqref{It:SP1}, \eqref{It:SP2} and \eqref{It:SP4} follow directly from the construction, while property \eqref{It:SP3} follows from the facts that $\area(\bigcap_n Y_n) = 0$. 
\end{proof} 

We are now ready to prove the main result.

\begin{proof}[Proof of Theorem~\ref{Thm:Rigidity}]

By Proposition~\ref{Prop:CombEquivalent}, there exists a pair of quasiconformally conjugate dynamically natural box mappings $F \colon \U \to \V$ and $\tilde F \colon \tilde \U \to \tilde \V$ such that $\Crit^\star(f) \subset \V$ and $\Crit^\star(\tf) \subset \tilde \V$. Put $H$ for this quasiconformal conjugacy, and let $K$ be its dilatation. Note that in this construction, the depth of $\V, \tV$ can be chosen arbitrarily large; call this the \emph{starting depth}.

Let $\Crit^{\#}(f)$ be the set of critical points in $\Crit^\star(f)$ that accumulate at at least one point in $\Crit^\star(f)$. By Theorem~\ref{Thm:Extract}, $\Crit^\#(f) \subset K(F)$ and $\Crit^\#(\tf) \subset K(\tilde F)$ (these are critical points that are `captured' by the box mappings, hence the notation with $\#$).

By construction, $\V$ and $\tV$ are nice neighborhoods of $\Crit^\#(f)$ and $\Crit^\#(\tf)$. By passing to puzzle pieces of $F$, $\tilde F$ of larger depth, we can construct a sequence of shrinking and nice (with respect to $f$) neighborhoods $\V_n$, $\tilde \V_n$ of these critical sets (they are shrinking by Theorem~\ref{Thm:BoxDynRigidity}). Here, $\V_1 = \V$, and the same for objects with the tilde. Restricting the quasiconformal conjugacy $H$ to these neighborhoods, we obtain a $K$-quasiconformal map $H \colon \V_n \to \tV_n$, $n \ge 1$, that respects the boundary marking induced by $\tau$.

By Theorem~\ref{Thm:BoxDynRigidity}, there exists a depth $\ft$ so that every critical point in $\Crit(f) \cap J_f$ is either in the puzzle graph of depth $\ft$, or it never lies on a puzzle graph and every two such critical points are separated by the puzzle graph of depth $\ft$. Moreover, the critical connected components of the complement of this graph have disjoint closures. Let us call $\ft$ the \emph{separation depth}. For $\tilde f$, the separation depth is also equal to $\ft$.

We can assume that the depths of puzzle pieces in $\V_n$ and $\tV_n$, $n \ge 1$, as puzzle pieces of $f$, are greater than the separation depth. By further increasing the starting depth, say to $\fm \ge \ft$, we can assume that the forward orbit of any point $c \in \Crit^{\star\sm\#}(f):= \Crit^\star(f) \sm \Crit^\#(f)$ does not intersect $\V_1$. 

For each pair $P, \tilde P$ of corresponding puzzle pieces of depth $\fm$, using Lemma~\ref{Lem:StartingQC} let us construct a $K'$-quasiconformal map $\phi_P \colon P \to \tilde P$, where $K'$ is uniform over all such pairs and each $\phi_P$ respects the boundary marking induced by the topological conjugacy $\tau$. This can be done because there are only finitely many puzzle pieces of depth $\fm$.

Now, let $\fm_n$ be the largest depth of puzzle pieces in $\V_n$. Let $\cQ_n$ be the puzzle neighborhood of $\Crit^{\star \sm \#}(f)$ of depth $\fm_n$. Pick the corresponding neighborhood $\tQ_n$ for $\tilde f$. For each connected component $Q$ of $\cQ_n$ and the corresponding component $\tilde Q$ of $\tQ_n$, let $f^k$, respectively $\tf^k$, be the iterate that maps $Q$, respectively $\tilde Q$ over a puzzle piece of depth $\fm$. The maps 
\[
f^k \colon Q \to f^k(Q), \quad \tf^k \colon \tilde Q \to \tf^k(\tilde Q)
\]      
are conformal isomorphisms by our choices of the depth. Hence, we can pull back the map $\phi_{f^k(Q)}$ to obtain a $K'$-quasiconformal map $\phi_{Q} \colon Q \to \tilde Q$. This map respects the boundary marking. Using such maps component-wise, we obtain a $K'$-quasiconformal homeomorphism $T_n \colon \cQ_n \to \tQ_n$ that respects the boundary marking. Note that $K'$ does not depend on the depth of $\cQ_n$. 

Combining $T_n$ and $H$, we obtain a $\max\{K,K'\}$-quasiconformal homeomorphism $G_n$ between nice, shrinking neighborhoods $\Y_n := \mathcal V_n \cup \mathcal Q_n$ of $\Crit^{\star}(f)$ and the corresponding neighborhoods $\tilde \Y_n$ of $\Crit^{\star}(\tf)$. These maps respect the boundary marking. By the Spreading Principle (Theorem~\ref{Thm:Spreading}), we can globalize this map to obtain a $\max\{K,K'\}$-quasiconformal map $\Phi_n \colon \Cc \to \Cc$ such that
\[
\tf \circ \Phi_n(z) = \Phi_n \circ f(z), \quad \forall z \in \Cc \,\sm\, \Y_n.
\] 
Let $\Phi$ be a sub-sequential limit of the sequence of uniformly quasiconformal maps $\Phi_n$. (Any two such subsequential limits must be equal as they coincide on an open and dense set of points.) Since the diameter of each component in $\Y_n$ shrinks to zero as $n \to \infty$, the map $\Phi$ provides the required quasiconformal conjugacy between $f$ and $\tf$. Since the measure of the set of points in $J_f$ that are never mapped to $\V_1$ is zero, the obtained map will be $K$-quasiconformal, with the support of this dilatation contained in $K(F)$. 

By Theorem~\ref{Thm:BoxLineFields}, since the box mapping $F$ is dynamically natural, $K(F)$ does not support any invariant line fields. Hence, $\Phi$ must be conformal. And since it is a global conformal map of $\Cc$ fixing $\infty$, it must be affine. Theorem~\ref{Thm:Rigidity} is proven. 
\end{proof}


\end{document}